\documentclass[preprint,3p,12pt]{elsarticle}
\usepackage{mathrsfs}
\usepackage{amsmath}
\usepackage{stmaryrd}
\usepackage{bbding}
\usepackage{dcolumn}
\usepackage{graphicx}
\usepackage{amsfonts}
\usepackage{amssymb}
\usepackage{psfrag}
\usepackage{wrapfig}
\usepackage{subfigure}
\usepackage{makeidx}
\usepackage{bm}
\usepackage{epsf}
\usepackage{epsfig}
\usepackage{setspace}
\usepackage{graphicx}
\usepackage{epstopdf}
\usepackage{psfrag}
\usepackage{subfigure}
\usepackage{color}

\begin{document}
\title{High-order Compact Gas-kinetic Scheme for Two-layer Shallow Water Equations on Unstructured Mesh }

\author[HKUST1]{Fengxiang Zhao}
\ead{fzhaoac@connect.ust.hk}

\author[HKUST1,HKUST2]{Jianping Gan}
\ead{magan@ust.hk}

\author[HKUST1,HKUST2,HKUST3]{Kun Xu\corref{cor}}
\ead{makxu@ust.hk}

\address[HKUST1]{Department of Mathematics, Hong Kong University of Science and Technology, Clear Water Bay, Kowloon, HongKong}
\address[HKUST2]{Center for Ocean Research in Hong Kong and Macau (CORE), Hong Kong University of Science and Technology, Clear Water Bay, Kowloon, Hong Kong}
\address[HKUST3]{Shenzhen Research Institute, Hong Kong University of Science and Technology, Shenzhen, China}
\cortext[cor]{Corresponding author}

\begin{abstract}
For the two-layer shallow water equations, a high-order compact gas-kinetic scheme (GKS) on triangular mesh is proposed.
The two-layer shallow water equations have complex source terms in comparison with the single layer equations.
 The main focus of this study is to construct a time-accurate evolution solution at a cell interface and to design a well-balanced scheme.
The evolution model at a cell interface provides not only the numerical fluxes, but also the flow variables.
The time-dependent flow variables at the closed cell interfaces can be used to update the cell-averaged gradients for the discretization of the  
the source terms inside each control volume in the development of the well-balanced scheme.
Based on the cell-averaged flow variable and their gradients, high-order initial data reconstruction can be achieved with compact stencils.
The compact high-order GKS has advantages to simulate the flow evolution in complex domain covered by unstructured mesh.
Many test cases are used to validate the accuracy and robustness of the scheme for the two-layer shallow water equations.
\end{abstract}

\begin{keyword}
Two-layer shallow water equations; Gas-kinetic scheme; High-order compact reconstruction; Unstructured mesh
\end{keyword}

\maketitle

\section{Introduction}

The shallow water equations (SWE) are useful in studying both large-scale ocean circulations and small-scale coastal and channel flows, such as tsunamis, pollutant transport, tidal waves, and dam break problems. However, real flows often exhibit stratification, which cannot be captured accurately by a single-layer SWE. For instance, the injection of freshwater into seawater creates plumes that are important for the coastal marine environment, with salinity stratification being a possible feature. In addition, the flow velocity at coastal area may vary significantly or exhibit stratification along the depth. To model the stratified water flow, the multi-layer SWE will be used
with the superposition of coupled layers with the force interaction between them.
This paper will focus on the development of high-order compact scheme for the
two-layer SWE (TLSWE), which is the basis for the multi-layer SWE.
In particular, the numerical scheme developed in this study for TLSWE can be naturally
extended to multi-layer SWE with the inclusion of the interaction between layers as the source term
and their dynamic effect in the calculation of numerical fluxes.

Numerous numerical schemes have been developed for solving SWE with second-order accuracy \cite{leveque-1998,zhou2001-surface,xu2002-swe}.
High-order numerical methods have gained popularity in recent years due to their advantages in accuracy and computational efficiency \cite{highorder-efficiency,wang2007high,wang2019-AIA,wang2021-AIA-benchmark}. As a result, several high-order numerical schemes have been proposed for solving SWE \cite{xing2006-DG,dambreak-2007,highorder-efficiency}.
However, there are few works on numerical methods for the two-layers SWE. Most of them are still based on the 1-D model \cite{castro2001-smallpertur,abgrall2009,bouchut2010-1dtlswe,2018nonflat-channel} or 2-D model on structured mesh \cite{kurganov2009-smallpertur,dudzinski2013}.
Second-order schemes for 2D TLSWE on unstructured mesh have been developed \cite{liu2021}, with great difficulties due to the loss of hyperbolicity under certain conditions and the stiff coupling between layers with the product of flow variables and their derivatives \cite{abgrall2009,castro2011-hyperbolicity,TLSWE-3forms,liu2021}.

Unstructured mesh is highly adaptable to complex geometries, making it a popular choice for numerical studies on real flow simulations \cite{SWE-DG}. This is particularly relevant for coastal hydrodynamics simulation, given the irregular and multiscale nature of coastal boundary geometries. However, to construct a high-order finite volume scheme on unstructured mesh presents a challenge due to the use of large stencils in the reconstruction \cite{dumbser}. Most high-order schemes for single-layer SWE on unstructured mesh are based on the discontinuous Galerkin (DG) formulation \cite{highorder-efficiency,SWE-DG}, and there are few high-order finite volume schemes for solving TLSWE.
The DG method updates the inner degrees of freedom (DOFs) from its weak formulation, and is widely used to solve compressible gas dynamics equations \cite{cockburn2,wangZJ2009CPR} due to its compact spatial discretization.
However, for flows with discontinuities, additional numerical treatments, such as identifying troubling cells and limiting procedures,
must be designed within the DG framework \cite{qiu2,qiu3,shu-review}.
In this study, a high-order compact gas-kinetic scheme (GKS) will be constructed.
The finite volume GKS updates both cell-averaged flow variables and their gradients from the moments of the time-accurate gas distribution function at a cell interface and compact initial reconstruction can be obtained.
At the same time, the multistage and multiderivative method will be adopted for achieving high-order temporal accuracy with less stages \cite{li2019AIA}.

The structure of this paper is as follows. Section 2 introduces the GKS for TLSWE.
Section 3 discusses the high-order compact reconstruction on unstructured mesh and temporal discretization.
In Section 4, the compact GKS is validated by studying shallow water flow in various cases.
Finally, Section 5 is the conclusion.

\section{Two-layer shallow water equations and gas-kinetic evolution model}

This section will present the gas-kinetic evolution model for solving TLSWE.
The corresponding GKS for TLSWE will be constructed based on the extension of the scheme for the single-layer SWE \cite{zhao2021-swe},
where the interaction between layers will be explicitly included in the scheme.

\subsection{Two-layer shallow water equations}

\begin{figure}[!htbp]
\begin{center}
\includegraphics[width=0.45\textwidth]{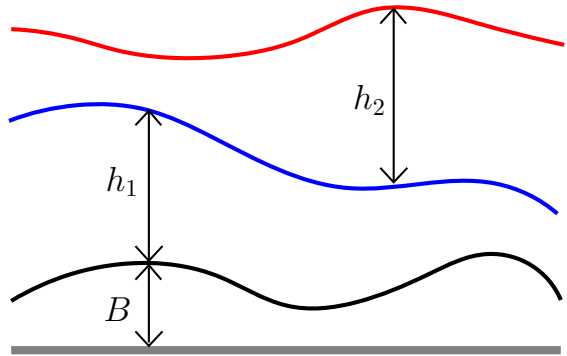}
\caption{Schematic of the two-layer shallow water flow.}\label{0-TLSWE-Schematic}
\end{center}
\end{figure}

In \cite{TLSWE-3forms}, three equivalent forms of TLSWE are presented.
In this study, the conservative form of TLSWE will be adopted and the interaction between layers is included in the source term,
\begin{equation}\label{SWE-macro}
\frac{\partial \textbf{W}}{\partial t}+ \frac{\partial \textbf{F}^x(\textbf{W})}{\partial x}+
\frac{\partial \textbf{F}^y(\textbf{W})}{\partial y}=\textbf{S}(\textbf{W}),
\end{equation}
where
\begin{equation*}
{\textbf{W}} =
\left(
\begin{array}{c}
h_2\\
h_2U_2\\
h_2V_2\\
h_1\\
h_1U_1\\
h_1V_1\\
\end{array}
\right), \\
{\textbf{F}^x} =
\left(
\begin{array}{c}
h_2U_2\\
h_2U_2^2+\frac{1}{2}Gh^2_2\\
h_2U_2V_2\\
h_1 {U_1}\\
h_1U_1^2+\frac{1}{2}Gh^2_1\\
h_1U_1V_1\\
\end{array}
\right),\\
{\textbf{F}^y} =
\left(
\begin{array}{c}
h_2V_2\\
h_2U_2V_2\\
h_2V_2^2+\frac{1}{2}Gh^2_2\\
h_1V_1\\
h_1U_1V_1\\
h_1V_1^2+\frac{1}{2}Gh^2_1\\
\end{array}
\right),
\end{equation*}
and
\begin{equation*}
{\textbf{S}} =
\left(
\begin{array}{c}
0\\
-Gh_2B_x-Gh_2h_{1,x}\\
-Gh_2B_y-Gh_2h_{1,y}\\
0\\
-Gh_1B_x-\chi Gh_1h_{2,x}\\
-Gh_1B_y-\chi Gh_1h_{2,y}\\
\end{array}
\right).
\end{equation*}
Here $\textbf{W}$ are the flow variables, and $\textbf{F}^x$ and $\textbf{F}^y$ are the corresponding fluxes in the $x$ and $y$ directions. $B$ is the bottom topography, $G$ is the gravitational acceleration, and $\chi$ is the density ratio defined as $\chi=\rho_2/\rho_1$, where $\rho_1$ and $\rho_2$ are the densities of the first and second fluid layer.
The flow variables of the lower and upper layers are denoted as $\mathbf{W}_1$ and $\mathbf{W}_2$, respectively. The fluxes of the two layers are $(\mathbf{F}_1^{x},\mathbf{F}_1^{y})$ and $(\mathbf{F}_2^{x},\mathbf{F}_2^{y})$ with the corresponding source terms $\mathbf{S}_1$ and $\mathbf{S}_2$. Fig.\ref{0-TLSWE-Schematic} presents a schematic of the two-layer shallow water flow.

By adopting the form of TLSWE in Eq. (\ref{SWE-macro}), the equations for each layer are similar as the single-layer SWE except
the additional source term related to the interaction between layers.
The source term makes the TLSWE conditionally hyperbolic \cite{TLSWE-3forms}, which may cause difficulty in the construction of numerical scheme based on Riemann solver and flux splitting method.
In addition, the source terms related to the interaction between layers are nonlinear, which have challenges in the
discretization for the high-order schemes.
In the gas-kinetic scheme, the dynamics in the TLSWE will be recovered by the time evolution of gas distribution function,
and the effect of source term will be incorporated into the particle transport process.
The numerical fluxes will be directly evaluated  from the time-dependent gas distribution function.
Since the governing equations of each layer in the TLSWE have similar forms, a general formulation for one of the layers will be presented in the following.

\subsection{Gas-kinetic evolution model}
The GKS is based on the time evolution solution of the gas distribution function for the flux evaluation \cite{xu2}.
The gas-kinetic BGK model can be written as \cite{xu2002-swe}
\begin{equation}\label{bgk-gravity}
f_t +\textbf{u}\cdot \nabla_{\textbf{x}} f +\nabla\Phi \cdot\nabla_{\textbf{u}} f=\frac{g-f}{\tau},
\end{equation}
where $f$ is the distribution function $f(\textbf{x},t,\textbf{u})$,  $\textbf{u}=(u,v)$ is the particle velocity,
and $g$ is the equilibrium state approached by $f$.
$\tau$ is the relaxation time.
$\nabla \Phi$ is the acceleration of particle due to external force and is related to the source term in TLSWE, such as the force from bottom topography and the friction.
The equilibrium state $g$ is a Maxwellian distribution function \cite{xu2002-swe},
\begin{equation}\label{g-SWE}
\begin{split}
g=h\big(\frac{\lambda}{\pi}\big)e^{-\lambda(\mathbf{u}-\mathbf{U})^2},
\end{split}
\end{equation}
where $\lambda$ is defined by $\lambda=1/Gh$. Due to the conservation in relaxation process from $f$ to $g$, $f$ and $g$ satisfy the compatibility condition,
\begin{equation}\label{bgk-compatibility}
\int \frac{g-f}{\tau}\pmb{\psi} \mathrm{d}\Xi=\textbf{0},
\end{equation}
where $\pmb{\psi}=(\psi_1,\psi_2,\psi_3)^T=(1,u,v)^T$ and $\text{d}\Xi=\text{d}u\text{d}v$.

Based on the moments of the gas distribution function, the flow variables and their fluxes can be obtained.
Due to the similar equations for different layers, the schemes for layer 1 and layer 2 can be formulated similarly.
In the general scheme, the macroscopic flow variables and the fluxes can be obtained from the distribution function $f$ as
\begin{equation}\label{bgk-g-to-W}
{\textbf{W}}
=\int f \pmb{\psi} \mathrm{d}\Xi,
\end{equation}
and
\begin{equation}\label{bgk-g-to-F}
{\big(\textbf{F}^x,\textbf{F}^y\big)^T}
=\int f \pmb{\psi} \mathbf{u} \mathrm{d}\Xi.
\end{equation}
The source term $\textbf{S}$ becomes
\begin{equation}\label{bgk-g-to-S}
{\textbf{S}}
=-\int \nabla\Phi\cdot\nabla_{\textbf{u}} f \pmb{\psi} \mathrm{d}\Xi,
\end{equation}
and $\nabla\Phi$ is determined by
\begin{equation}\label{bgk-phi}
\nabla\Phi=\mathbf{S}/h,
\end{equation}
where $\mathbf{S}$ takes $\mathbf{S}_1=h_1(0,-GB_x-G\chi h_{2,x},-GB_y-G\chi h_{2,y})^\mathrm{T}$ and $\mathbf{S}_2=h_2(0,-GB_x-Gh_{1,x},-GB_y-Gh_{1,y})^\mathrm{T}$ for layer 1 and layer 2, respectively.

The formal solution of the BGK model in Eq. (\ref{bgk-gravity}) with external forcing term is
\begin{equation}\label{f-integral-solution}
f(\textbf{x},t,\textbf{u})=\frac{1}{\tau}\int_0^t g(\textbf{x}^{'},t',\textbf{u}^{'})e^{-(t-t')/\tau}\mathrm{d}t'
+e^{-t/\tau}f_0(\textbf{x}_0,\textbf{u}_0),
\end{equation}
where $\textbf{x}$ is the numerical quadrature point on the cell interface for flux evaluation, and $\textbf{x}$ can be set as $(0,0)$ for simplicity in a local coordinate system with both normal and tangential directions as the x- and y-directions.
The formal solution describes an evolution process for the distribution function. The trajectory of fluid particle is given by  $\textbf{x}=\textbf{x}^{'}+\textbf{u}^{'}(t-t^{'})+\frac{1}{2}\nabla\Phi(t-t^{'})^2$, and the velocity of the particle is $\textbf{u}=\textbf{u}^{'}+\nabla\Phi(t-t^{'})$.
The acceleration has a second-order effect ($\sim t^2$) on the particle trajectory, but
has the first-order contribution ($\sim t$) to the particle velocity.

The second-order in time and the ell-balanced explicit evolution solution $f$ is obtained  for SWE \cite{zhao2021-swe}.
In this paper, the same evolution solution of $f$ is used for the individual layer.
The solution of $f$ is
\begin{equation}\label{SWE-2nd-order}
\begin{split}
f(\textbf{x},t,\textbf{u})&=\overline{g}(\mathbf{x},0,\mathbf{u})\big[ C_1+ C_2 \big( \overline{\mathbf{a}}^l \cdot\mathbf{u}H(u) +\overline{\mathbf{a}}^r \cdot\mathbf{u}(1-H(u)) \big) +C_3\overline{A} \big]\\
                            &+C_2\overline{g}(\mathbf{x},0,\mathbf{u}) \big[-2 \alpha_{k,m} \overline{\lambda} \big( \nabla\Phi^l H(u)+\nabla\Phi^r(1-H(u)) \big) \cdot (\mathbf{u}-\overline{\mathbf{U}}) \big]\\
                            &+C_4\big[g^l(\mathbf{x},0,\mathbf{u})H(u)+g^r(\mathbf{x},0,\mathbf{u})(1-H(u))\big] \\
                            &+C_5g^l(\mathbf{x},0,\mathbf{u})\big[\mathbf{a}^l\cdot \mathbf{u} -2 \alpha_{k,m} \lambda^l \nabla\Phi^l \cdot(\mathbf{u}-\mathbf{U}^l) \big]H(u) \\
                            &+C_5g^r(\mathbf{x},0,\mathbf{u})\big[\mathbf{a}^r\cdot \mathbf{u} -2 \alpha_{k,m} \lambda^r \nabla\Phi^r \cdot(\mathbf{u}-\mathbf{U}^r) \big](1-H(u)),
\end{split}
\end{equation}
where $\alpha_{k,m}~(k=1,2,~m=1,2,3)$ are constants for a well-balanced scheme,
$(\alpha_{1,1},\alpha_{1,2},\alpha_{1,3})=(1,3/4,1/4)$ and $\alpha_{2,m}=1$, $m$ and $k$ is related to taking moment, and the details are given in the Appendix of \cite{zhao2021-swe}.
The coefficients $C_i~(i=1,2,\cdots,5)$ are
\begin{equation*}
\begin{split}
C_1&=1-e^{-t/\tau}, ~ C_2=-\tau(1-e^{-t/\tau})+te^{-t/\tau}, ~ C_3=-\tau(1-e^{-t/\tau})+t, \\
C_4&=e^{-t/\tau}, ~ C_5=-t e^{-t/\tau}.
\end{split}
\end{equation*}
The fluxes at the cell interface are evaluated by taking moments of the above gas distribution function
and the total transport of mass and momentum within a time step can be further integrated in time.
More details in the formulation can be found in \cite{xu2}.

\subsection{Acceleration force modeling at the interface between two water layers}
The interaction between layers is modeled as the acceleration term in the kinetic equation.
The spatial derivatives of the water column height determine the acceleration,
where the values of the height derivatives can be obtained by the compact reconstruction at the cell interface.
However, the possible discontinuity of the interface can trigger a sudden ``pull'' or ``push'' between water layers.
For cases with discontinuities, the spatial derivatives of the water height from the reconstruction will not be used to calculate the force,
and the ``step effect'' due to the discontinuity needs to be considered.
The acceleration from a discontinuous interface will be modeled.

\begin{figure}[!htb]
\centering
\includegraphics[width=0.35\textwidth]{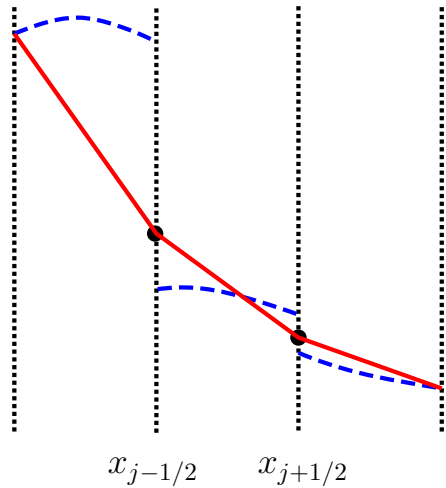}
\caption{\label{modeling-discontinuity-Schematic} A 1-D schematic for the modeling discontinuous interface between water layers to evaluate the acceleration term in the BGK model. In the 1-D discretized space, for the reconstructed distributions of the water interface inside each cell with the possible discontinuities at the cell interface,
the topography for evaluating the source term between layers is represented by the connected solid line. }
\end{figure}

Without loss of generality, for the momentum equation of layer 2 as an example,
the corresponding acceleration is given by
\begin{align*}
\nabla \Phi_2=-GB_x -\chi Gh_{1,x}.
\end{align*}
For the continuous bottom topography $B$ and water height $h_1$,
the acceleration can be directly evaluated based on the functions of $B_x$ and the reconstructed $h_{1,x}$.
However, when the water height is discontinuous at a cell interface, such as the reconstructed dash lines in Fig. \ref{modeling-discontinuity-Schematic}, the corresponding forcing term between layers will be modeled from
a re-constructed continuous profile at the cell interface.
The construction of this continuous profile will take into account the forcing interaction between neighboring cells.

Firstly, let's construct the continuous line at the cell interface.
In each cell, the continuous line connects the respective unique values of the water height on the cell interfaces $x_{j\pm1/2}$, which are denoted by $\widehat{h}_1(x_{j\pm1/2})$ as the black dots in Fig. \ref{modeling-discontinuity-Schematic} with the values given later.
The continuous line in the cell is obtained as
\begin{equation*}
P_j^1(x)=\frac{1}{2}(\widehat{h}_1(x_{j+1/2})+\widehat{h}_1(x_{j-1/2}))+\frac{\widehat{h}_1(x_{j+1/2})-\widehat{h}_1(x_{j-1/2})}{x_{j+1/2}-x_{j-1/2}}(x-x_{j}),
\end{equation*}
where $P_j^1(x)$ is a linear interpolation based on the values at the cell interfaces of the cell.
$\widehat{h}_1(x_{j+1/2})$ are modeled based on the discontinuous left and right states 
\begin{equation*}
\widehat{h}_1(x_{j+1/2})=\xi h^l_1(x_{j+1/2}) +(1-\xi)h^r_1(x_{j+1/2}),
\end{equation*}
where $h^{l}_1(x_{j+1/2})$ and $h^r_1(x_{j+1/2})$ are the reconstructed values at the cell interface.
$\xi$ is a coefficient for the convex combination, and it is defined as
\begin{equation*}
\xi=\frac{1}{2}\mathrm{erfc}\big( (U^{l}_1(x_{j+1/2})+U^{r}_1(x_{j+1/2}))/2 \big),
\end{equation*}
where the function $\mathrm{erfc}(\cdots)$ is the complementary error function,
and $U^{l,r}_1(x_{j+1/2})$ are the left and right values of the velocity at the interface.
 $\mathrm{erfc}(\cdots)$ makes a smooth transition from $2$ to $0$ when the independent variable covers $(-\infty,+\infty)$ with a value $\mathrm{erfc}(0)=1$.
The above linear distribution in the cell has dynamically upwind-biased slope.

In the smooth case, the updated derivative of the water height can be used in the evaluation of acceleration inside each cell.
In order cope with both discontinuous and smooth cases, the final derivative of the water height is determined 
by the following nonlinear convex combination method
\begin{align}
\begin{split}
\widetilde{h}^l_{1,x}(x_{j+1/2})&=w_{j+1/2} h^l_{1,x}(x_{j+1/2}) +(1-w_{j+1/2})P^1_{j,x}(x_{j+1/2}),\\
\widetilde{h}^r_{1,x}(x_{j+1/2})&=w_{j+1/2} h^r_{1,x}(x_{j+1/2}) +(1-w_{j+1/2})P^1_{j+1,x}(x_{j+1/2}),
\end{split}
\end{align}
where $w_{j+1/2}$ is a nonlinear weighting function to identify the smoothness of the solution.
$w_{j+1/2}$ tends to $1$ in the smooth region and to $0$ in the discontinuous region. 
The value of $w_{j+1/2}$ is the same nonlinear weight as that in the high-order time stepping reconstruction scheme of \cite{zhao2023direct}. 

In two dimensions, similar modeling of the derivative of the water height can be done.
Different from the one-dimensional one, the modeled continuous line in Fig. \ref{modeling-discontinuity-Schematic} is extended to a 2-D continuous plane.
A smooth linear interpolation in the cell is determined by the following constraints.
\begin{equation*}
P^1(x^c_{k},y^c_{k})=\widehat{h}_1(x^c_{k},y^c_{k}), ~k=1,2,3,
\end{equation*}
where $(x^c_{k},y^c_{k})$ is the center of the cell interface, $\widehat{h}_1(x^c_{k},y^c_{k})$ can be obtained
by taking the arithmetic average of the values $\widehat{h}_1(x_{m},y_{m})$, where $m=1,2$, on the Gaussian quadrature points of the corresponding cell side.

%

\section{Compact GKS based on high-order compact reconstruction}
In this section, the compact GKS for the TLSWE will be constructed, where the high-order compact reconstruction to obtain the initial values of flow distributions is implemented and the two-stage fourth-order (S2O4) temporal discretization is used.
Since  the two layers in the shallow water equations can be numerically treated in the same way, the evolutions for $\mathbf{W}_1$ and $\mathbf{W}_2$ will be presented by the discretization of $\mathbf{W}$ below.

\subsection{Finite volume discretization}
Taking moments $\pmb{\psi}$ on Eq. (\ref{bgk-gravity}), the flow variables in a cell $\Omega_j$  are updated by
\begin{equation}\label{SWE-semi-discrete}
\frac{\partial \textbf{W}_{j}}{\partial t}=-\frac{1}{\big|\Omega_j\big|} \int_{\partial \Omega_j} \textbf{F}\cdot \textbf{n} \mathrm{d}l +\frac{1}{\big|\Omega_j\big|} \int_{\Omega_j} \textbf{S} \mathrm{d} \Omega_j,
\end{equation}
where $\textbf{W}_{j}$ is the cell-averaged flow variable, $\textbf{F}=(\textbf{F}^x,\textbf{F}^y)$ is the time-dependent flux at cell interface,
which can be obtained from the moments of the gas distribution function in Eq. (\ref{SWE-2nd-order}).
The $\textbf{W}_{j}$ is defined as
\begin{align}\label{SWE-cell-average}
\textbf{W}_{j}&\equiv \frac{1}{\big| \Omega_j \big|} \int_{\Omega_j} \textbf{W}(\textbf{x}) \text{d} \Omega.
\end{align}
The line integral of the flux in Eq. (\ref{SWE-semi-discrete}) can be discretized by a q-point Gaussian quadrature formula,
\begin{align}\label{SWE-semi-flux}
\begin{split}
-\frac{1}{\big|\Omega_j\big|} \int_{\partial \Omega_j} \textbf{F}\cdot \textbf{n} \mathrm{d} l &= -\frac{1}{\big|\Omega_j\big|} \sum_{l=1}^{l_0}\big( \big|\Gamma_{l} \big| \sum _{k=1}^q \omega_k \textbf{F}(\textbf{x}_k)\cdot \textbf{n}_l \big) \\
& \equiv \mathcal{L}^F_j(\textbf{W}),
\end{split}
\end{align}
where $\big|\Gamma_{l} \big|$ is the side length of the cell, $l_0$ is the total number of cell sides, such as $l_0=3$ for a triangular mesh, $\textbf{n}_l$ is the unit outer normal vector, and $q$ and $\omega_k$ are the total number of integration points and weights of the Gaussian integration formula.
In order to evaluate the above numerical flux, the initial data $\textbf{W}(\textbf{x}_k)$ is reconstructed using the compact spatial stencil, which are presented in Section 3.3.
The cell-averaged $\textbf{S}$ becomes
\begin{equation}\label{SWE-semi-source}
\begin{split}
\frac{1}{\big|\Omega_j\big|} \int_{\Omega_j} \textbf{S} \mathrm{d} \Omega_j \equiv \mathcal{L}^S_j(\textbf{W}).
\end{split}
\end{equation}

\subsection{Discretization for source term}
The source term in the momentum equations includes two parts, the first one depends on the bottom topography,
and the second one is related to the variation of the interface between layers and the water height of the up layer.

The first part of the source term depending on the bottom topography is defined as
\begin{equation}\label{SWE-semi-source}
\begin{split}
\frac{1}{\big|\Omega_j\big|} \iint_{\Omega_j} \textbf{S} \mathrm{d} \Omega_j &= h_j(0,-GB_{j,x},-GB_{j,y})^T \\
& \equiv \mathcal{L}^{S_1}_j(\textbf{W}),
\end{split}
\end{equation}
where $h_j$ is the cell average of $h$ in $\Omega_j$.
High-order spatial and temporal discretizations of the first part can be implemented directly, as in the single-layer SWE in \cite{zhao2021-swe}.

The second part of $\mathcal{L}^{S}_j(\textbf{W})$ is related to the variation of the water height.
Taking the source term in the equation of $h_1U_1$ as an example, the spatial discretization becomes
\begin{align}\label{SWE-source-discretization}
\begin{split}
\mathcal{L}^{S_2}_j(\textbf{W}) &\equiv \frac{1}{|\Omega_j|}\int_{\Omega_j} -\chi Gh_1h_{2,x} \mathrm{d}x \mathrm{d}y \\
&=-\chi G \frac{1}{|\Omega_j|}\int_{\Omega_j} h_1 \mathrm{d}x \mathrm{d}y \cdot \frac{1}{|\Omega_j|}\int_{\Omega_j} h_{2,x} \mathrm{d}x \mathrm{d}y +O(\Delta X^2) \\
&=-\chi G \frac{1}{|\Omega_j|}\int_{\Omega_j} h_1 \mathrm{d}x \mathrm{d}y \cdot \frac{1}{|\Omega_j|}\int_{\partial \Omega_j} h_{2}n_x \mathrm{d}\Gamma +O(\Delta X^2) \\
&=-\chi G \sum_{k=1}^{3}\sum_{l=1}^{2} \widetilde{w}_{k,l}h_1(\mathbf{x}_{k,l}) \cdot \frac{1}{|\Omega_j|} \sum_{k=1}^{3} \big(\sum_{l=1}^{2} w_{k,l} h_{2}(\mathbf{x}_{k,l}) \big)n_{k,x}|\Gamma_k| +O(\Delta X^2).
\end{split}
\end{align}
where $|\Omega_j|$, $\big|\Gamma_{k} \big|$, $\omega_l$, $n_x$ and $\mathbf{x}_{k,l}$ have the same definition as those in Eq. (\ref{SWE-semi-flux}),
$\Delta X$ is the mesh cell size, $\widetilde{w}_{k,l}$ is the weight to obtain the numerical integration over $\Omega_j$
based on $h_1(\mathbf{x}_{k,l})$, and $\widetilde{w}_{k,l}=1/6$.
The second-order spatial discretizations is implemented in Eq. (\ref{SWE-source-discretization}).
High-order discretization of $\mathcal{L}^{S_2}_j(\textbf{W})$ can be achieved by introducing more numerical integration points.
However, considering the balance between accuracy and efficiency, the simple method given in Eq. (\ref{SWE-source-discretization})
is adopted for the spatial discretization of the second part of the source term in this paper.

The compact GKS of the TLSWE is a well-balanced scheme.
The well-balanced property is achieved through the balance of the  time-accurate flux function at the cell interface
and the spatial discretization of the source terms inside the control volume.
In the previous study \cite{zhao2022compact}, the well-balanced GKS for the single-layer SWE on triangular mesh has been developed, where
a corresponding well-balanced evolution solution of the gas distribution function shown in Eq. (\ref{SWE-2nd-order}) is obtained.
For the TLSWE, with the well-balanced initial conditions
\begin{align*}
\begin{split}
&h_1+B=\mathrm{Const}, \\
&h_2=\mathrm{Const},
\end{split}
\end{align*}
and $(U_1,V_1)=(U_2,V_2)=(0,0)$,
at the quadrature points on the cell interface, the initial conditions should be
$\nabla h_2=\mathbf{0}$ and $\nabla (h_1+B)=\mathbf{0}$.
With the adoption of water level reconstruction technique \cite{zhou2001-surface}, this initial condition can be preserved numerically.
As a result, the compact GKS for the TLSWE can keep such a solution and the scheme is a well-balanced one.
In the following, the solution update in the compact GKS on the triangular mesh will be presented.

\subsection{The time evolutions of flow variables and their derivatives}
By adopting the S2O4 time stepping method \cite{li2019AIA,pan1},
the fully discretized form of the TLSWE in Eq. (\ref{SWE-semi-source}) over the cell $\Omega_j$ in a time step $[t^n,t^{n+1}]$ is given by
\begin{equation}\label{SWE-fully-discrete}
\begin{split}
\textbf{W}^{n+1/2}_j=&\textbf{W}^n_j+\frac{1}{2}\Delta t\mathcal
{L}_j(\textbf{W}^n)+\frac{1}{8}\Delta t^2\frac{\partial}{\partial
	t}\mathcal{L}_j(\textbf{W}^n),\\
\textbf{W}^{n+1}_j=&\textbf{W}^n_j+\Delta t\mathcal{L}_j(\textbf{W}^n)+\frac{1}{6}\Delta t^2\frac{\partial}{\partial t}\mathcal{L}_j(\textbf{W}^n)
+\frac{1}{3}\Delta t^2\frac{\partial}{\partial t}\mathcal{L}_j(\textbf{W}^{n+1/2}),
\end{split}
\end{equation}
where $\mathcal{L}_j=\mathcal{L}^F_j+\mathcal{L}^S_j$ includes the flux and source term contribution.

\begin{figure}[!htb]
\centering
\includegraphics[width=0.40\textwidth]{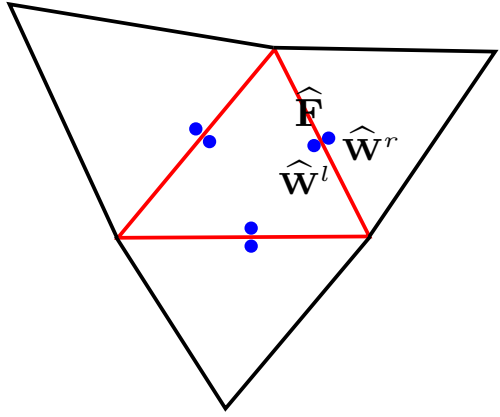}
\caption{\label{0-compact-update} {A schematic of flow variables and fluxes on the interface given by the time-dependent evolution solution in the compact GKS. The dots represent time accurate flow variables at the interface with possible discontinuities, but corresponding a
single valued flux function.  }}
\end{figure}

In the current compact GKS, besides the update of cell-averaged flow variables in Eq. (\ref{SWE-fully-discrete}),
the cell-averaged derivatives can be updated as well by the Gauss's theorem as
\begin{equation}\label{slope}
\nabla \mathbf{W}_j(t^{n+1}) =\frac{1}{|\Omega_j|} \int_{\partial \Omega_j} {\bf W}(\mathbf{x},t^{n+1}) {\bf n} \mathrm{d} S,
\end{equation}
with the discretized form
\begin{align}\label{slope-1}
\nabla \textbf{W}_{j}^{n+1}& =\frac{1}{\big| \Omega_j \big|}\sum_{l=1}^{l_0} \big(|\Gamma_l| \mathbf{n}_{l} \sum _{k=1}^q \omega_k\textbf{W}^{n+1}(\textbf{x}_k) \big),
\end{align}
where $|\Omega_j|$, $\big|\Gamma_{l} \big|$, $l_0$, $\omega_k$ and $\textbf{n}_l$ have the same definition as those in Eq. (\ref{SWE-semi-flux}). The flow variables $\mathbf{W}(\mathbf{x},t^{n+1})$ should be provided at the inner sides of the cell boundary of the control volume at the time step $t^{n+1}$.
Fig.\ref{0-compact-update} shows the time-accurate flow variables and fluxes on the cell interface from the evolution solution of the gas distribution function in the compact GKS.
In the discrete scheme, the discontinuous evolution solution $\mathbf{W}^{l,r}(\mathbf{x},t^{n+1})$ at the cell interface have been obtained
in the GKS for the highly compressible Navier-Stokes solutions \cite{zhao2023direct}.
However, in the current study for the shallow water equations, a continuous evolution solution, namely $\mathbf{W}^{l}=\mathbf{W}^{r}$, for the update the cell-averaged derivatives within the cell by Eq. (\ref{slope}) seems work very well.
In order to obtain a high-order time-accurate flow variable at the quadrature point in Eq. (\ref{slope-1}),
the macroscopic flow variable is evolved by two stages
\begin{equation}\label{S2O3}
\begin{split}
\textbf{W}^{n+1/2}(\textbf{x})=\textbf{W}^n(\textbf{x})+\frac{1}{2}\Delta t \textbf{W}_{t}^n(\textbf{x}), \\
\textbf{W}^{n+1}(\textbf{x})  =\textbf{W}^n(\textbf{x})+\Delta t \textbf{W}_{t}^{n+1/2}(\textbf{x}).
\end{split}
\end{equation}

\subsection{High-order compact reconstruction}
\begin{figure}[!htb]
\centering
\includegraphics[width=0.40\textwidth]{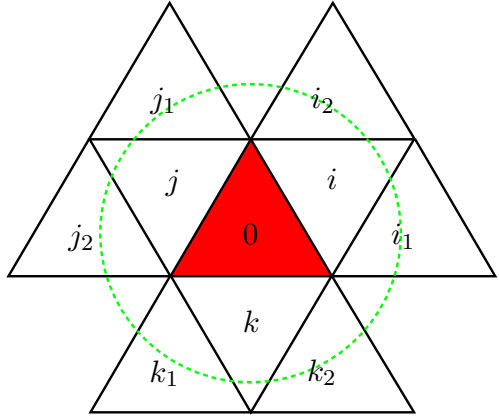}
\caption{\label{0-stencil-2d} A schematic of reconstruction stencil of compact GKS. The  dotted circle is the physical domain of dependence
of the cell 0, which is covered compactly by the cells surrounding 0. In each cell, the flow variables and their derivatives in the x- and y-directions are known.}
\end{figure}

In this section, the high-order compact spatial reconstruction for flow variables will be presented.
Based on the cell averages and their derivatives, the high-order reconstruction can be obtained compactly with the  stencils involving the closest neighboring cells only, as shown in Fig. \ref{0-stencil-2d}.
The compact stencil provides consistent domains of dependence between the numerical and physical ones.
The reconstruction with the accuracy from fourth-order to sixth-order can be obtained on the compact stencils \cite{zhao2022compact}.
The fourth-order reconstruction will be used in this study.

For the fourth-order reconstruction, $P^3(\bm{x})$ polynomial is constructed as
\begin{align}\label{Recons-poly}
\begin{split}
P^3(\bm{x})=\sum_{k=0}^{9} a_k \varphi_k(\bm{x}),
\end{split}
\end{align}
where $a_k$ is the degrees of freedom (DOFs) of $P^3(\bm{x})$,  the total number of $a_k$ is $10$ and the complete polynomial basis with the highest order of $3$ are included, and $\bm{x}=(x,y)$ is the coordinate.
The basis function $\varphi_k(\bm{x})$ can take the zero-averaged basis as
\begin{align}\label{Recons-eqs}
1, ~\delta x-\overline{\delta x}^{(0)}, ~\delta y-\overline{\delta y}^{(0)}, ~\frac{1}{2}\delta x^2-\overline{\frac{1}{2}\delta x^2}^{(0)}, ~\delta x\delta y-\overline{\delta x\delta y}^{(0)}, ~\frac{1}{2}\delta y^2-\overline{\frac{1}{2}\delta y^2}^{(0)}, ~\cdots.
\end{align}

To fully determine $P^3(\bm{x})$, the DOFs on the cells of the compact stencil is selected to give the constraints on $P^3(\bm{x})$.
\begin{align}\label{Recons-eqs}
\begin{split}
&\big(\frac{1}{\big|\Omega_l \big|} \int_{ \Omega_l} \varphi_k(\bm{x}) \mathrm{d}x \mathrm{d}y \big) a_k=Q_l, \\
&\big(\frac{1}{\big|\Omega_l \big|} \int_{ \Omega_l} \varphi_{k,x}(\bm{x}) \mathrm{d}x \mathrm{d}y \big) a_k=Q_{l,x}, \\
&\big(\frac{1}{\big|\Omega_l \big|} \int_{ \Omega_l} \varphi_{k,y}(\bm{x}) \mathrm{d}x \mathrm{d}y \big) a_k=Q_{l,y},
\end{split}
\end{align}
where the same subscript $k$ of $\varphi_k$ and $a_k$ on the left-hand side of the equations follow the Einstein summation.
$Q_{l}$, $Q_{l,x}$ and $Q_{l,y}$ are the DOFs in the cells for any component of $\mathbf{W}$.

Due to arbitrary geometrical triangular mesh, the number of equations $M$ in Eq. (\ref{Recons-eqs}) should be greater than the number of DOFs $a_k$ to avoid an ill-conditioned system. For the fourth-order reconstruction, the set of DOFs $S_0$ is given by
\begin{align}\label{Recons-stencil}
S_0=\{Q_{l_1},Q_{l_2,x},Q_{l_2,y}\}, ~l_1=0,i,j,k,i_1,i_2,\cdots,k_2,~l_2=0,i,j,k.
\end{align}
Eq. (\ref{Recons-eqs}) can determine a linear system of $a_k$, and it is written as
\begin{equation} \label{Recons-system-0}
a_0=Q_0,
\end{equation}
and
\begin{equation} \label{Recons-system}
\left(
\begin{array}{cccc}
A_{1,1} & A_{1,2} & \cdots & A_{1,9} \\
A_{2,1} & A_{2,2} & \cdots & A_{2,9} \\
\vdots  &\vdots   & \vdots & \vdots  \\
A_{9,1} & A_{9,2} & \cdots & A_{9,9} \\
\widetilde{A}^x_{0,1} & \widetilde{A}^x_{0,2} & \cdots & \widetilde{A}^x_{0,9} \\
\widetilde{A}^y_{0,1} & \widetilde{A}^y_{0,2} & \cdots & \widetilde{A}^y_{0,9} \\
\vdots  &\vdots   & \vdots & \vdots  \\
\widetilde{A}^y_{3,1} & \widetilde{A}^y_{3,2} & \cdots & \widetilde{A}^y_{3,9} \\
\end{array}
\right)
\left(
\begin{array}{c}
a_1 \\
a_2 \\
\vdots \\
a_9 \\
\end{array}
\right)
=
\left(
\begin{array}{c}
Q_{1}-Q_{0} \\
Q_{2}-Q_{0} \\
\vdots      \\
Q_{9}-Q_{0} \\
Q_{0,x}h    \\
Q_{0,y}h    \\
\vdots      \\
Q_{3,y}h    \\
\end{array}
\right) ,
\end{equation}
where $A_{l,k}$ and $\widetilde{A}_{l,k}$ are defined as
\begin{align}
\begin{split}
A_{l,k}&=\frac{1}{\big|\Omega_l \big|} \int_{ \Omega_l} \varphi_k(\bm{x}) ~\mathrm{d}x \mathrm{d}y,\\
\widetilde{A}^i_{l,k}&=\frac{h}{\big|\Omega_l \big|} \int_{ \Omega_l} \partial \varphi_{k}(\bm{x})/\partial r ~\mathrm{d}x \mathrm{d}y, ~r=x,y, ~k=1,2,\cdots,9.\\
\end{split}
\end{align}
The system can be solved by the least square (LS) method. The solution of $a_k,~(k=1,2,\cdots,9)$ is given by
\begin{align}
\mathbf{a}=\big[(\mathbf{A}^\mathrm{T} \mathbf{A})^{-1}\mathbf{A}^\mathrm{T}\big] \mathbf{Q},
\end{align}
where $\mathbf{a}$ is the vector of DOFs without $a_0$, $\mathbf{A}$ is the coefficient matrix in Eq. (\ref{Recons-system}), and $\mathbf{Q}$ is the vector of the RHS in Eq. (\ref{Recons-system}).

To deal with discontinuities in the solution, the nonlinear reconstruction is needed. The nonlinear compact reconstruction is obtained based on the WENO method by nonlinearly combining the high-order polynomial $P^3$ and several lower-order polynomials, where the lower-order polynomials are determined based on the sub-stencils by using the LS method. The nonlinear reconstruction in the compact GKS has been developed in \cite{zhao2022compact}, and the same techniques will be used here.

\section{Numerical validations}
The compact GKS for the two-layer SWE will be validated by the cases of two-layer shallow flow in this section. All the computations in this section are performed on 2-D triangular mesh.
The time step used in the computation is determined by the CFL condition as $\Delta t=CFL\frac{\Delta X}{U_{max}}$, where $\Delta X$ is the size of the mesh cell, $U_{max}=\mathrm{max}\{\sqrt{U_1^2+V_1^2}+\sqrt{Gh_1},\sqrt{U_2^2+V_2^2}+\sqrt{Gh_2}\}$, and $CFL$ number takes $0.5$.
The gravitational acceleration is taken as $G=9.81$ if not specified.

The collision time $\tau$ in the BGK model for inviscid flow at a cell interface is defined by
\begin{align*}
\tau=\varepsilon \Delta t + \varepsilon_{num}\displaystyle|\frac{h_l^2-h_r^2}{h_l^2+h_r^2}|\Delta t,
\end{align*}
where $\varepsilon=0.05$, $\varepsilon_{num}=5$, and $h^2_l$ and $h^2_r$ are the pressures at the left and right sides of a cell interface.
The reason for including the pressure jump term in the relaxation time is to enhance the artificial dissipation in case of bore wave.

\subsection{Accuracy test}
The accuracy of the compact GKS with high-order compact reconstruction is tested.
In order to calculate the error in the numerical solution, an initial condition with analytical evolution solution is used directly
\begin{align*}
\begin{split}
& h_1=0.9+0.02e^{-50((x-1)^2+(y-1)^2)}, \\
& h_2=1-h_1,
\end{split}
\end{align*}
with a uniform velocity $(U_1,V_1)=(U_2,V_2)=(1,1)$.
The density ratio is taken as $\chi=1.0$. The gravitational acceleration is $G=9.81$. The free boundary condition is taken.
The analytical solution of this problem is given by
\begin{align*}
& h_1(t)=0.9+0.02e^{-50((x-1-t)^2+(y-1-t)^2)}, \\
& h_2(t)=1-h_1(t), \\
& (U_1(t),V_1(t))=(U_2(t),V_2(t))=(1,1).
\end{align*}
The computational domain is taken as $[0,2]\times[0,2]$. The triangular mesh is used.

The $L^1$ errors of $h_1$ and $h_2$ at $t=0.1$ and the convergence orders are presented in Table \ref{Accu-test}.
The convergence order of current compact GKS does not keep the 4th order, which is due to the second-order approximation is used when discretizing the source term in Eq. (\ref{SWE-source-discretization}).
Although the optimal 4th-order convergence is not realized, the advantages of high resolution from the compact spatial reconstruction will be demonstrated in other complex flow problems.

\begin{table}[!h]
	\small
	\begin{center}
		\def\temptablewidth{0.7\textwidth}
		{\rule{\temptablewidth}{0.40pt}}
        \footnotesize
		\begin{tabular*}{\temptablewidth}{@{\extracolsep{\fill}}c|cc|cc}
         $h_{re}$ & $Error_{L^1}(h_1)$ &$\mathcal{O}_{L^1}(h_1)$  &$Error_{L^1}(h_2)$ &$\mathcal{O}_{L^1}(h_2)$    \\
			\hline
            1/8   & 1.2664e-04    & ~     & 1.2015e-04  & ~      \\
            1/16  & 1.9987e-05    & 2.66  & 1.8203e-05  & 2.72   \\
            1/32  & 1.0744e-06    & 4.22  & 9.4509e-07  & 4.27   \\
            1/64  & 1.4263e-07    & 2.91  & 1.2522e-07  & 2.92   \\
		\end{tabular*}
		{\rule{\temptablewidth}{0.40pt}}
	\end{center}
	\vspace{-6mm} \caption{\label{Accu-test} Accuracy test: errors and convergence orders of compact GKS with high-order reconstruction. }
\end{table}

\subsection{Well-balanced property}
The well-balanced property of the compact GKS on unstructured mesh is validated in the following.
The initial condition is a two-dimensional steady state solution with non-flat bottom topography.
The bottom topography is
\begin{align*}
& B(x,y)=0.5e^{-50[(x-1)^2+(y-1)^2]}.
\end{align*}
The steady state is
\begin{align*}
\begin{split}
& h_1=0.8-B(x,y), \\
& h_2=0.2,
\end{split}
\end{align*}
and all the velocities are $0$. The density ratio and the gravitational acceleration are taken as $\chi=1.0$ and $G=9.81$, respectively.
The computational domain is $[0,2]\times[0,2]$. The triangular mesh with cell size $\Delta X=0.05$ is used. The wall boundary condition is imposed on all the boundaries.

\begin{figure}[!htb]
\centering
\includegraphics[width=0.45\textwidth]{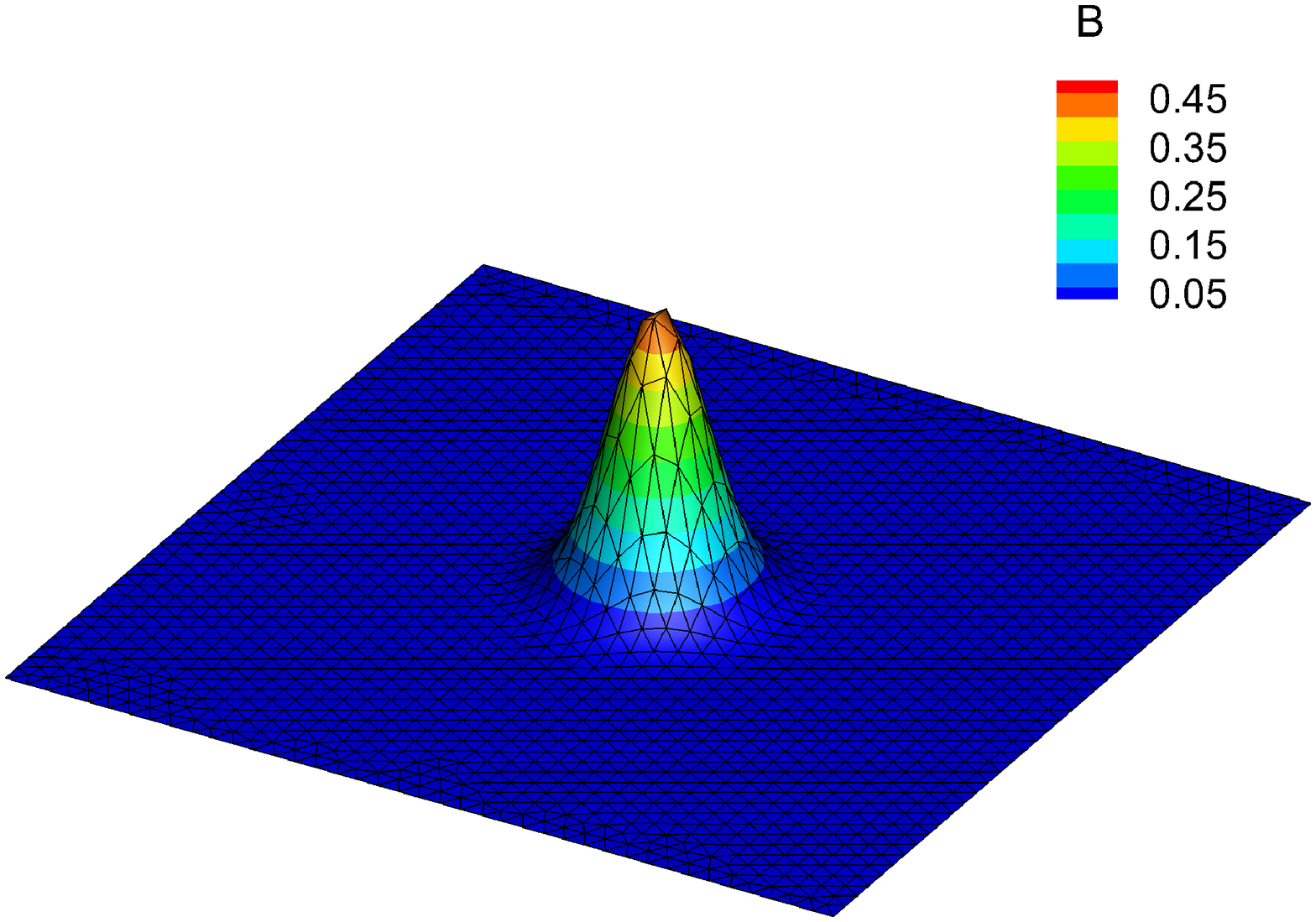}
\includegraphics[width=0.45\textwidth]{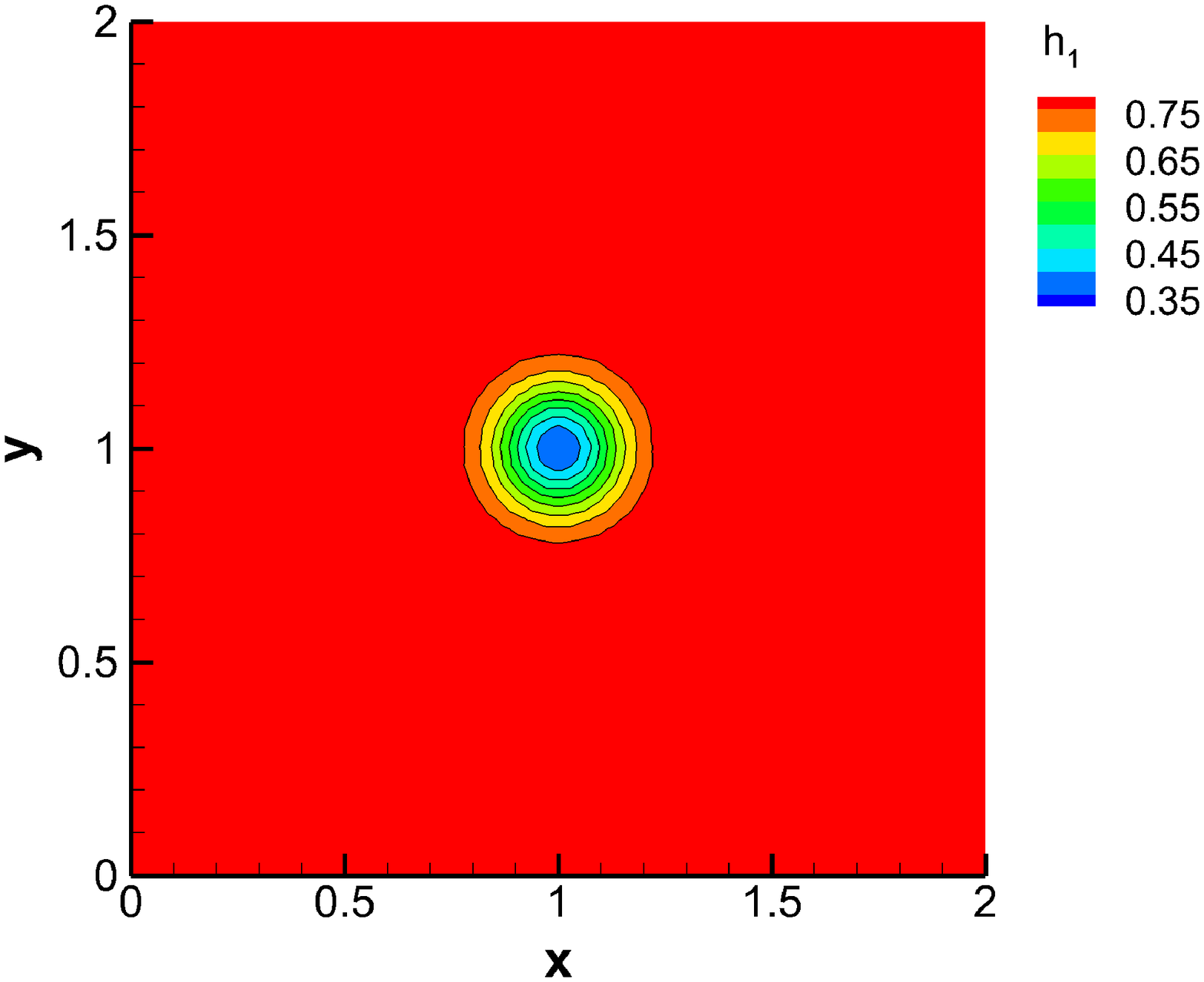}
\caption{\label{2d-well-balance-1} Well-balanced property study: the left figure shows the bottom profile and the unstructured mesh with cell size $\Delta X=0.05$, and the right figure presents the water level contours at $t=100$.}
\end{figure}

\begin{figure}[!htb]
\centering
\includegraphics[width=0.45\textwidth]{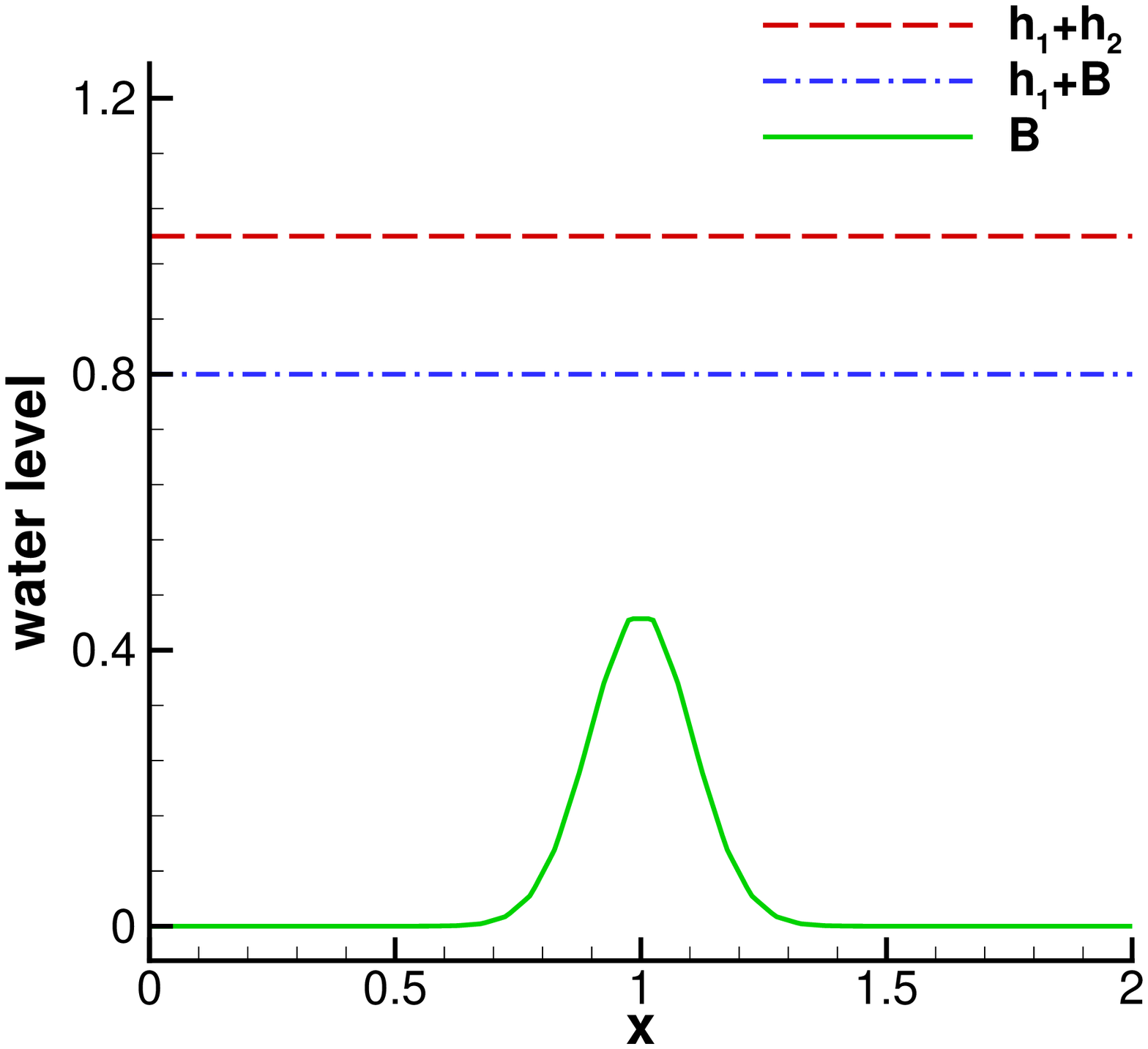}
\includegraphics[width=0.45\textwidth]{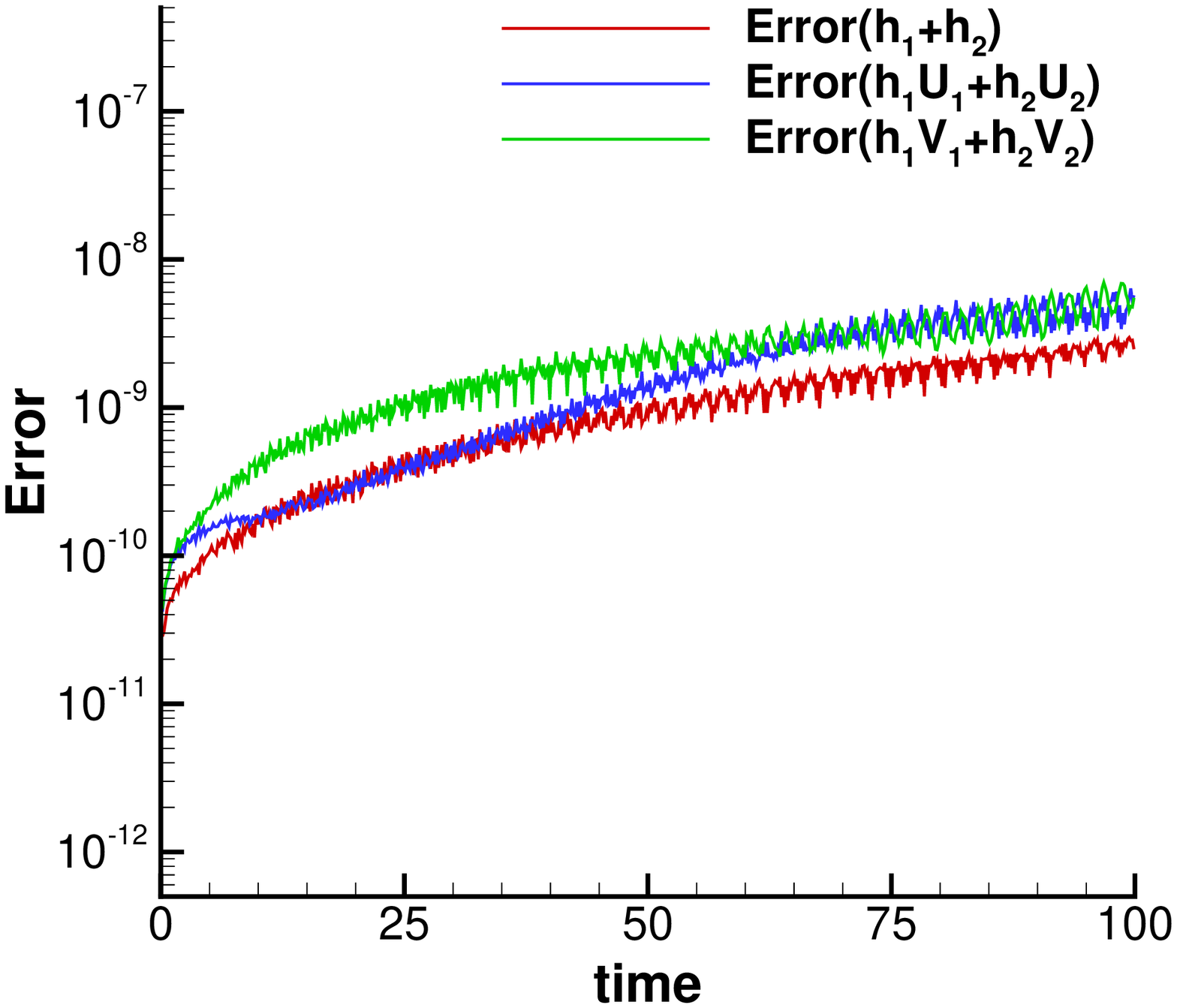}
\caption{\label{2d-well-balance-2} Well-balanced property study: the left figure shows the water level distributions along the horizontal centerline in the computational domain at $t=100$ and the right figures presents the time evolution of the numerical error.}
\end{figure}

The discretized bottom topography is shown in Fig. \ref{2d-well-balance-1}. The errors history of flow variables is plotted in Fig. \ref{2d-well-balance-2}. The error remains at the same level at different computational time.
At very long computation times, the errors of water surface level and momentum are less than $1.0\times10^{-8}$. The current compact GKS is able to maintain an initial balanced steady state solution.

\begin{figure}[!htb]
\centering
\includegraphics[width=0.45\textwidth]{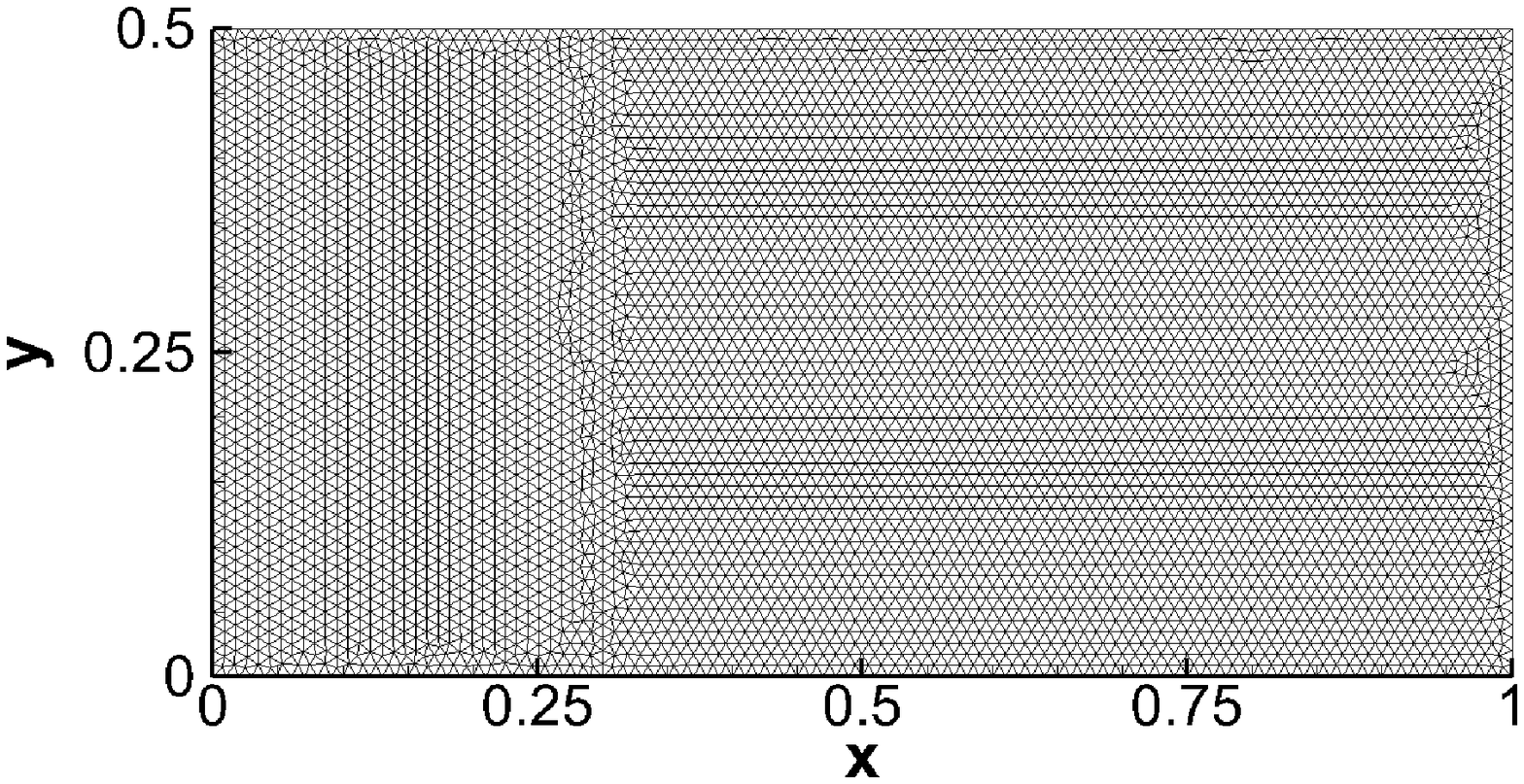}
\includegraphics[width=0.45\textwidth]{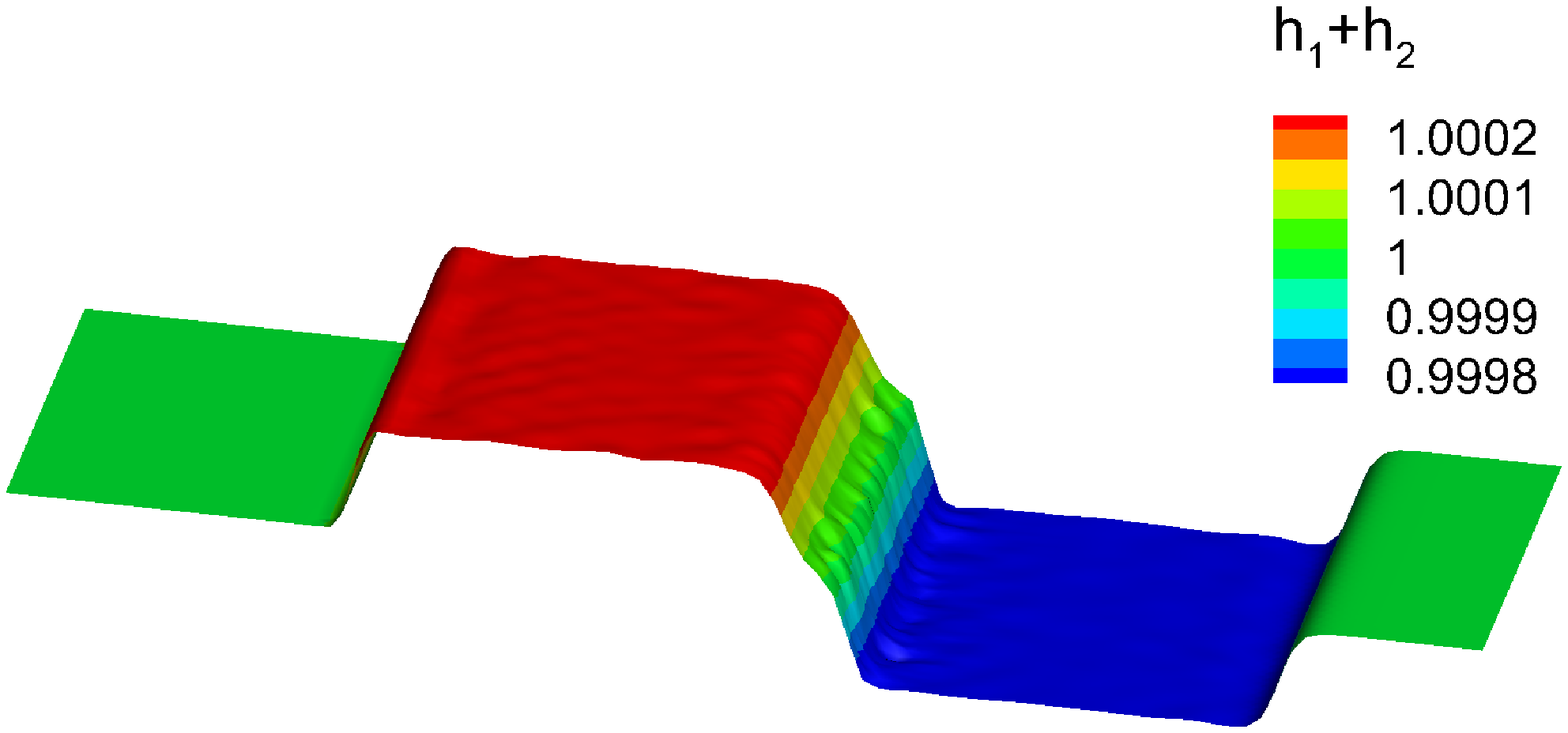}
\caption{\label{weak-dis-1} Riemann problem I: the left figure shows the unstructured mesh with cell size $\Delta X=1/100$,
and the right figure is the 3-D water surface distributions at $t=0.1$.}
\end{figure}

\begin{figure}[!htb]
\centering
\includegraphics[width=0.45\textwidth]{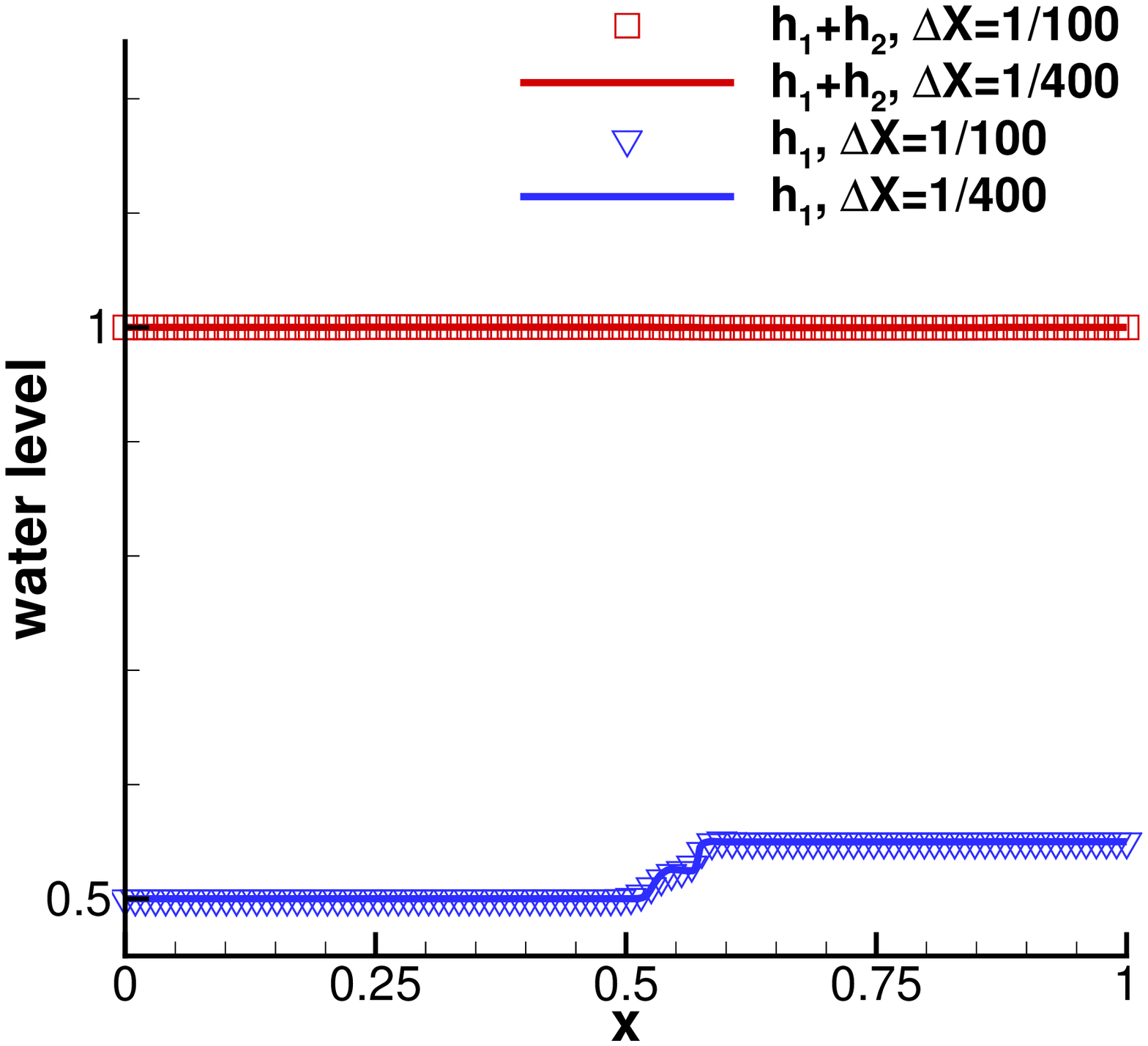}
\includegraphics[width=0.45\textwidth]{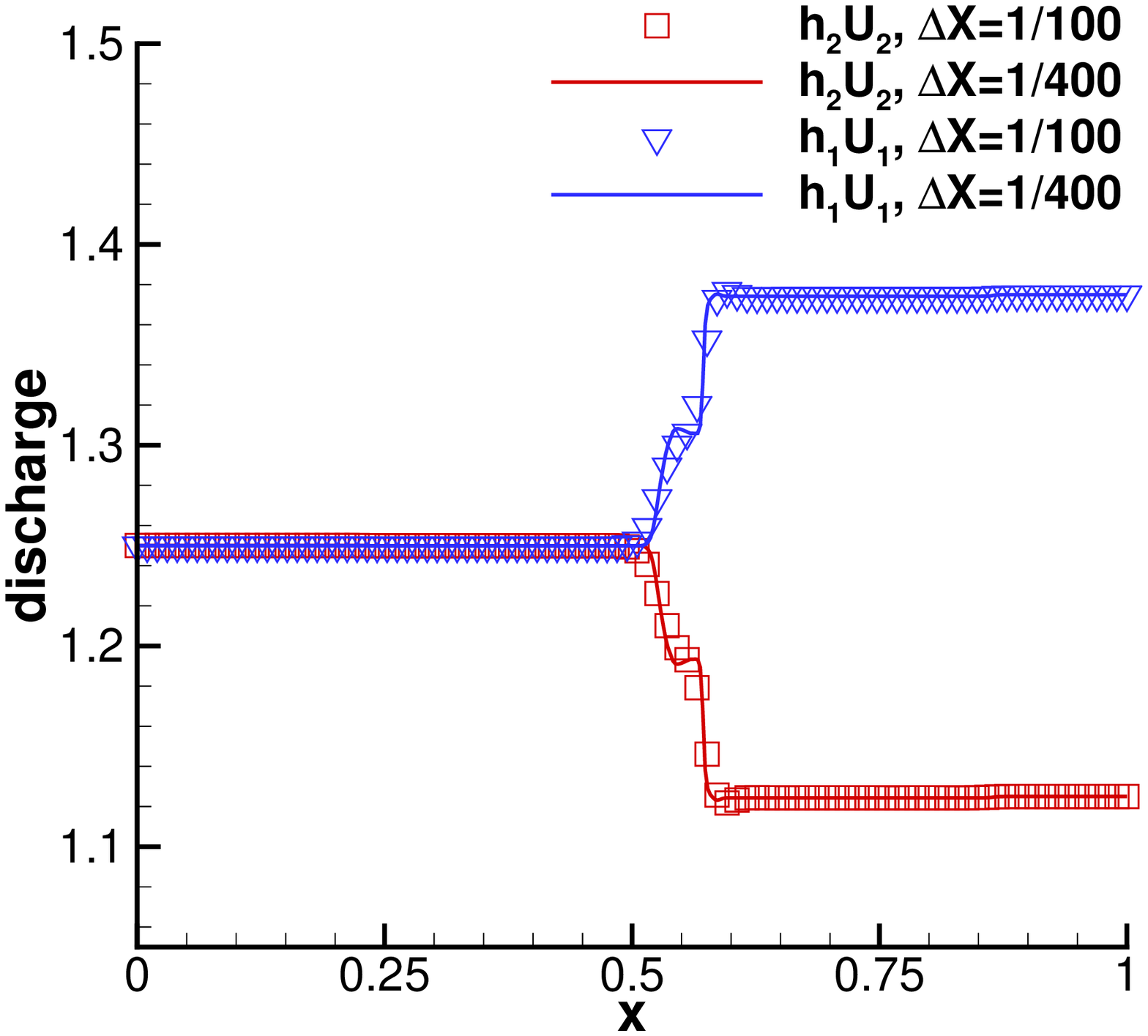}
\caption{\label{weak-dis-2} Riemann problem I: the left figure is the 1-D distributions of water level and the right figure is about the
discharge along the horizontal centerline at $t=0.1$.}
\end{figure}

\begin{figure}[!htb]
\centering
\includegraphics[width=0.45\textwidth]{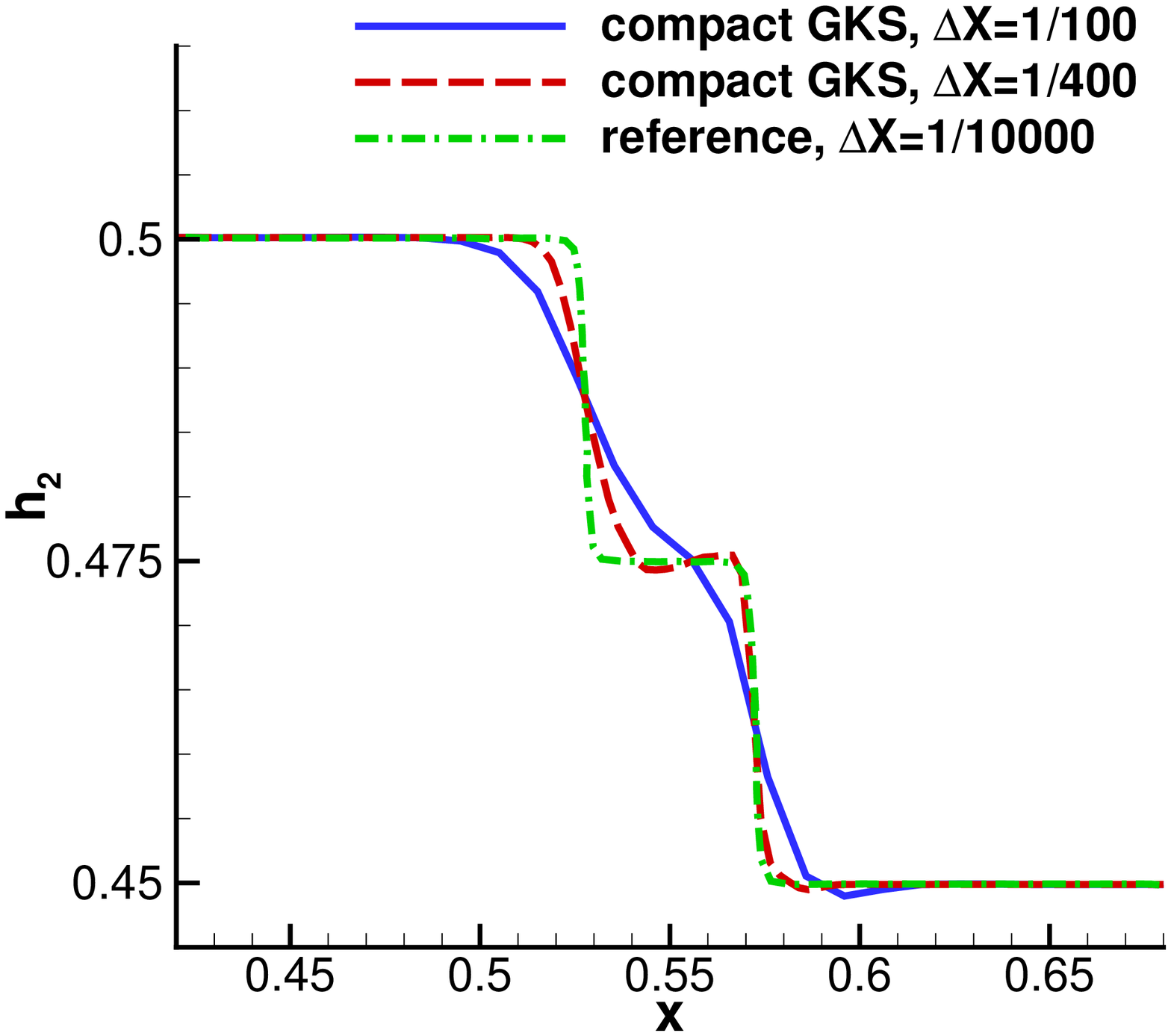}
\includegraphics[width=0.45\textwidth]{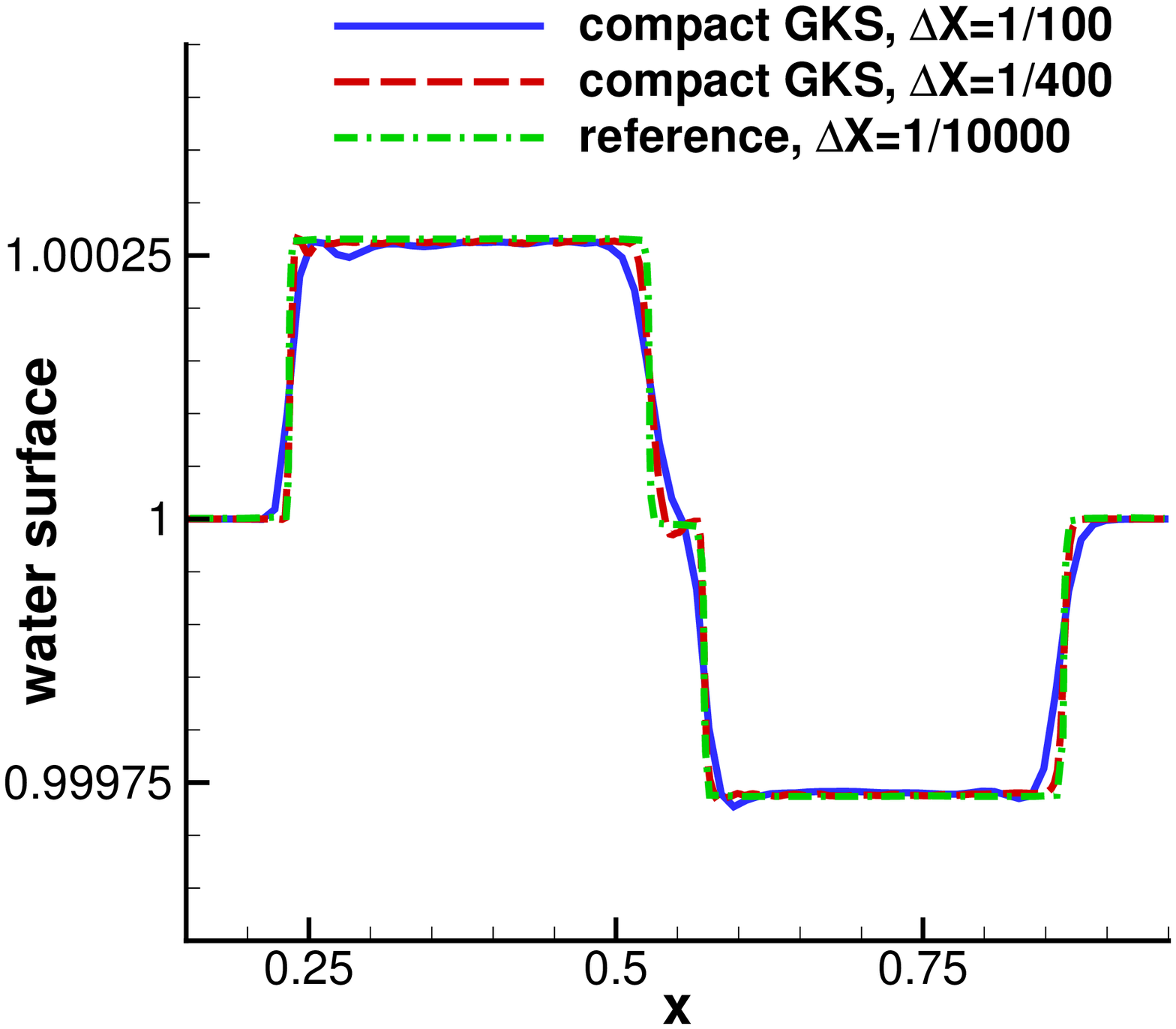}
\caption{\label{weak-dis-3} Riemann problem I: comparison of the water levels obtained by the compact GKS and from the reference solution.}
\end{figure}

\subsection{Riemann problems of TLSWE}

In this section, the Riemann problems with a discontinuity at the interface between two fluid layers are studied to validate the compact GKS for TLSWE. Due to unequal densities of the two layers of fluid, the discontinuity at the interface will evolve and propagate.

The first test was introduced to verify the stability of the numerical schemes for unsteady two-layer exchange flows \cite{castro2001-smallpertur}. It can also be used to evaluate the accuracy of different numerical schemes in computing unsteady solutions over a flat bottom.
The initial water level is set as
\begin{equation*}
(h_1,h_2) = \begin{cases}
(0.5,0.5),~~~~ 0\leq x<0.3,\\
(0.55,0.45), 0.3\leq x\leq1,
\end{cases}
\end{equation*}
and the uniform velocity $(U_1,V_1)=(U_2,V_2)=(2.5,0)$ is given in the whole domain.
In the computation, the 2-D computational domain is taken as $[0,1]\times[0,0.5]$, and the triangular mesh is used.
The computational time is $t=0.1$. The density ratio is $\chi=0.98$. The gravitational acceleration is taken as $G=10$ in this case.

The coarse mesh with $\Delta X=1/100$ used in the computation and the 3-D water surface obtained by the compact GKS are shown in Fig. \ref{weak-dis-1}.
The solution of the evolved free surface has a square-wave structure with small variation.
The current compact GKS captures this solution with no obvious numerical oscillations.
In Fig. \ref{weak-dis-2} and Fig. \ref{weak-dis-3}, the water levels along the horizontal centerline of the computational domain is plotted, where the results on a finer mesh with $\Delta X=1/400$ are also given to verify the mesh convergence solution from the current compact GKS.
To quantitatively verify the correctness of the results obtained by the current scheme, the reference solution obtained by the 1-D model with a cell size of $\Delta X=1/10000$ in \cite{kurganov2009-smallpertur} is also plotted.
The compact GKS gives consistent solutions on both coarse and fine meshes.
The resolution of the local solution structure on the fine mesh by the compact GKS is comparable to the reference solution.

\begin{figure}[!htb]
\centering
\includegraphics[width=0.45\textwidth]{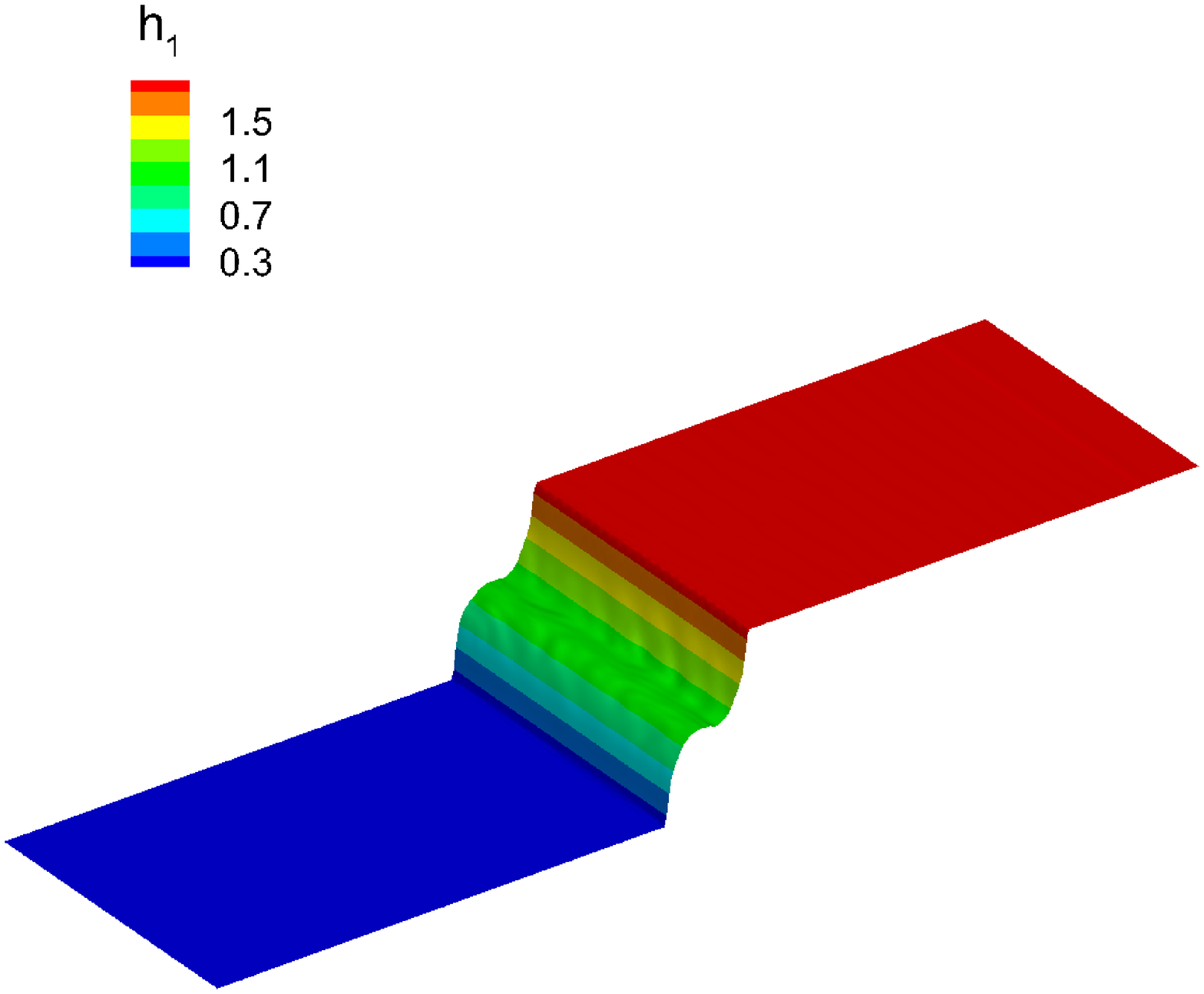}
\includegraphics[width=0.45\textwidth]{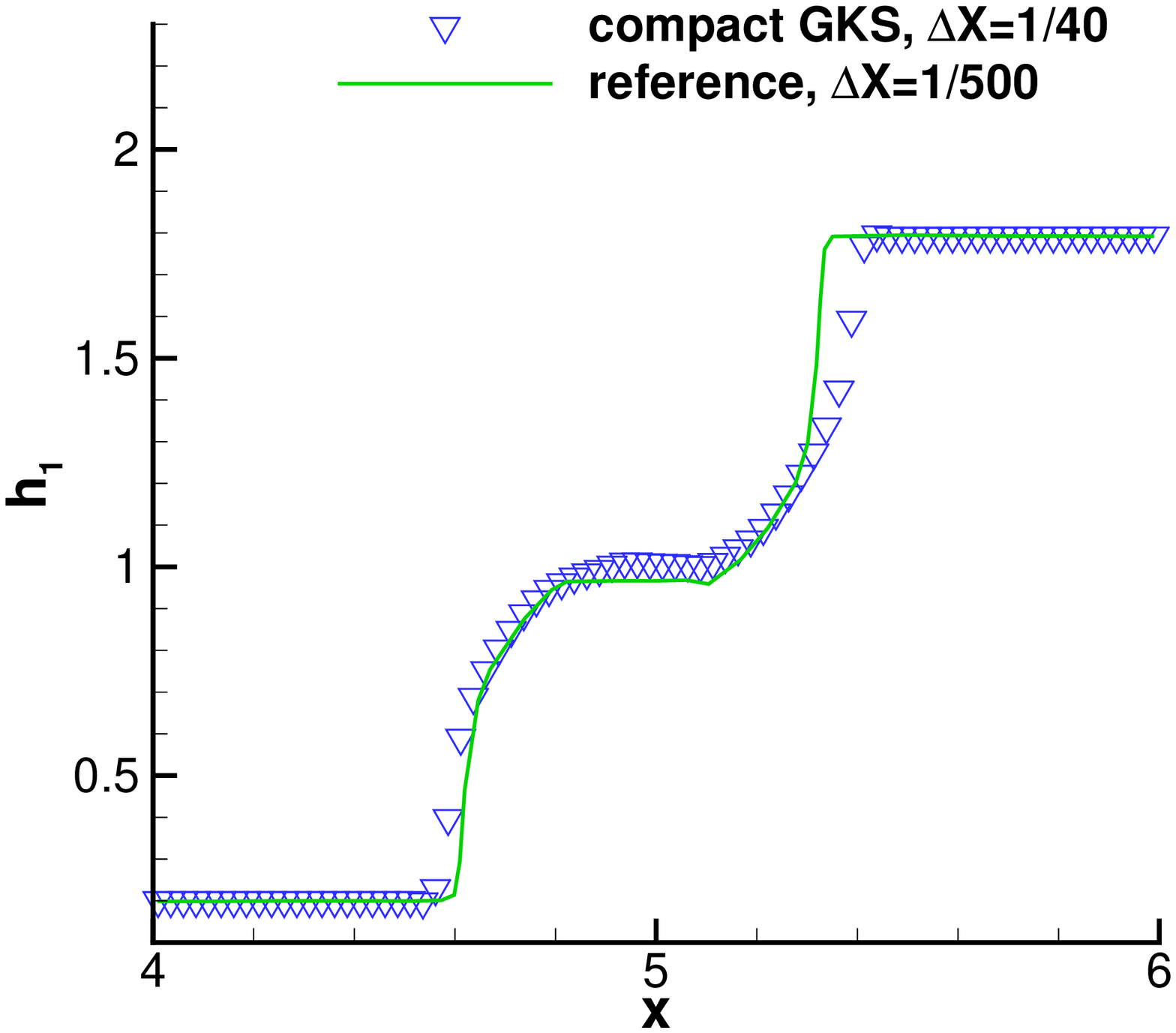}\\
\caption{\label{str-dis-1} Riemann problem II: the left figure shows the 3-D water level contours of $h_1$ and the right figure gives the 1-D distribution of $h_1$ along horizontal centerline at $t=1$.}
\end{figure}

\begin{figure}[!htb]
\centering
\includegraphics[width=0.45\textwidth]{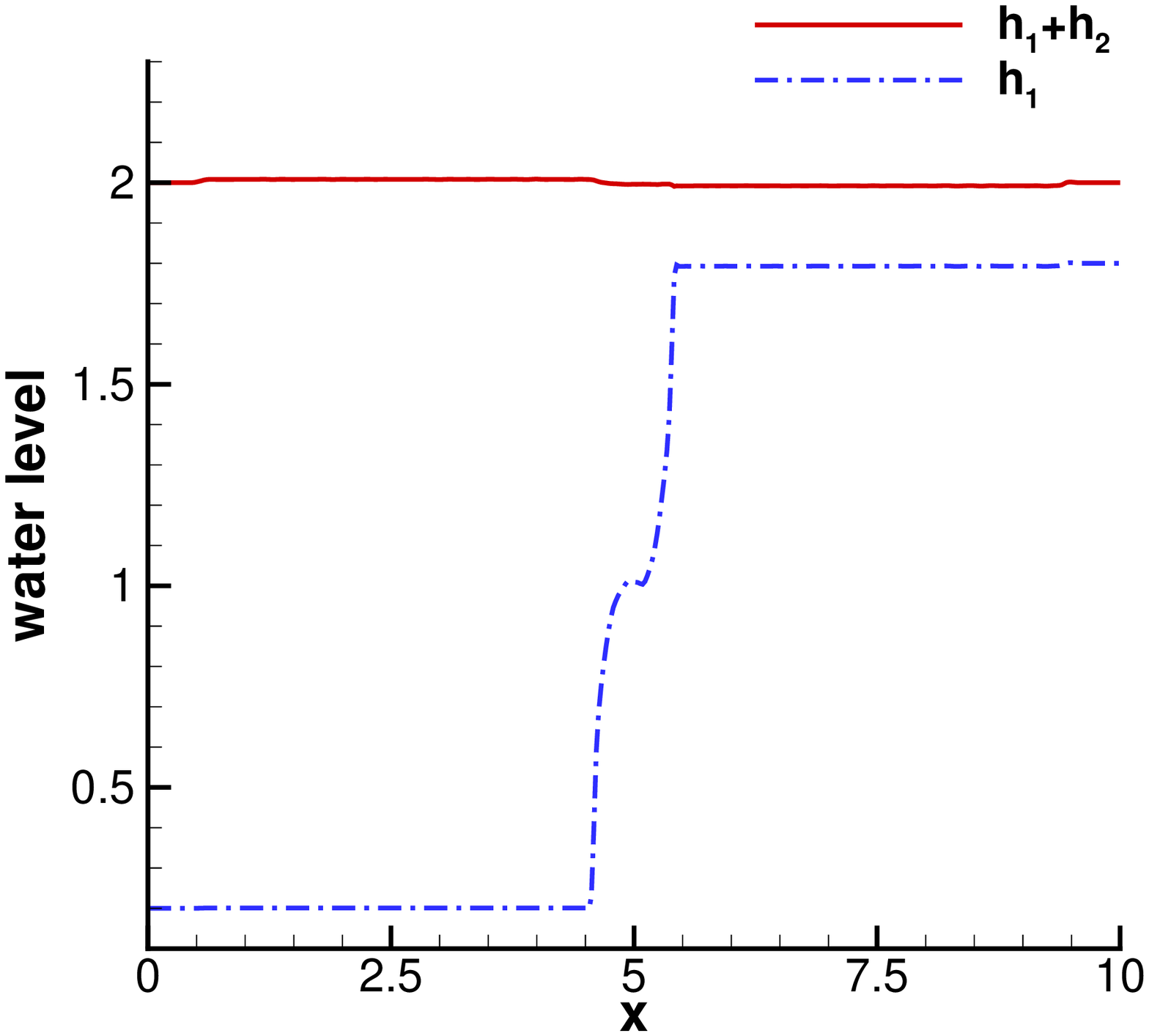}
\includegraphics[width=0.45\textwidth]{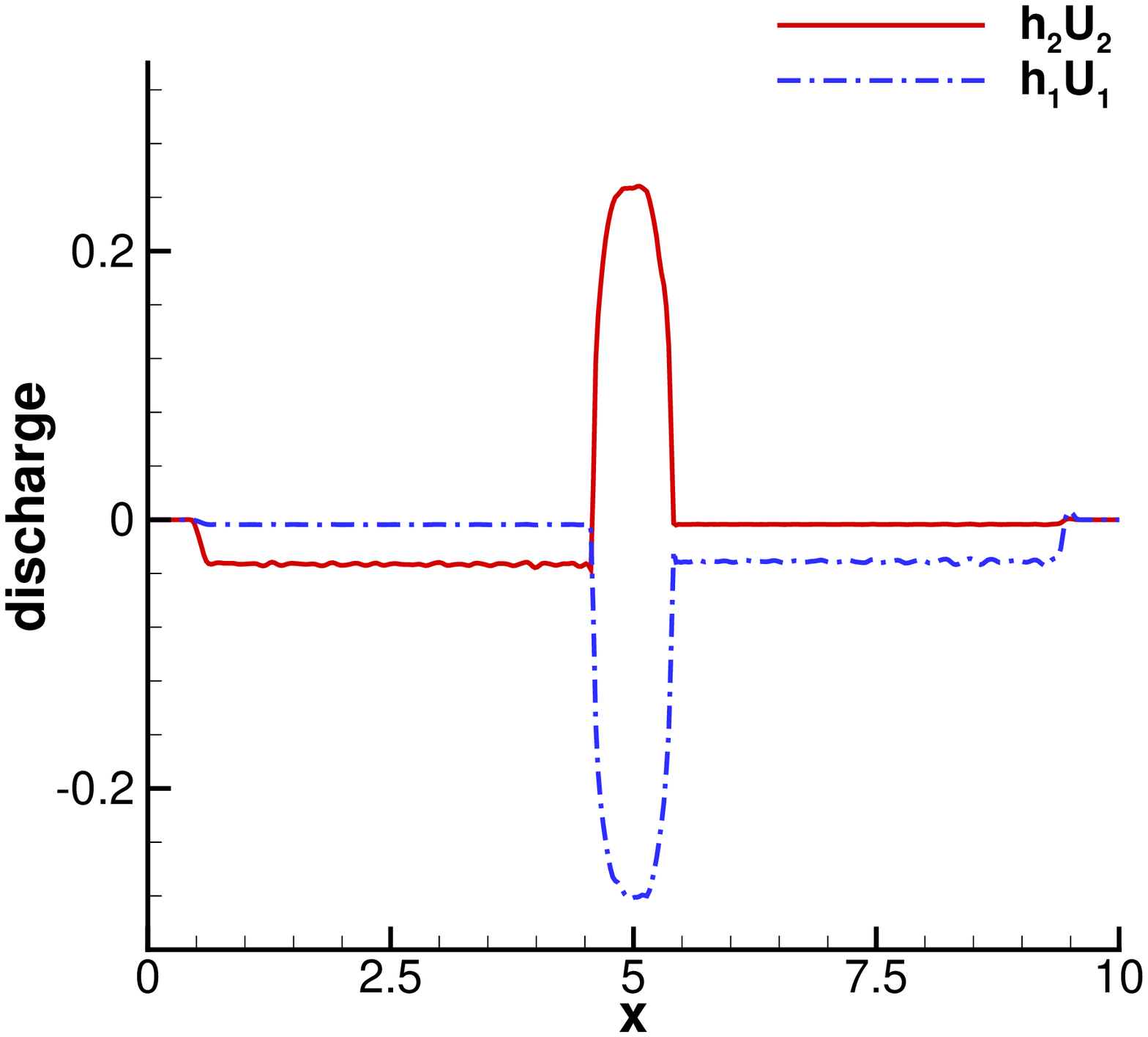}
\caption{\label{str-dis-2} Riemann problem II: the left figure is the 1-D distributions of water level and the right figure is about the
discharge along the horizontal centerline at $t=1$.}
\end{figure}

The second case is the Riemann problem with a large discontinuity at the interface between the two layers \cite{kurganov2009-smallpertur}.
The initial value of water levels is given by
\begin{equation*}
(h_1,h_2) = \begin{cases}
(0.2,1.8),  ~0 \leq x<5,\\
(1.8,0.2),  ~5 \leq x<10.
\end{cases}
\end{equation*}
The initial velocity is $0$. The water density ratio is $\chi=0.98$.
The gravitational acceleration is taken as $G=9.81$. The computational domain is set as $[0,10]\times[0,1]$. The triangular mesh with a cell size of $\Delta X=1/40$ is used in the computation.

The evolved results at $t=1.0$ obtained by the compact GKS is presented in Fig. \ref{str-dis-1} and Fig. \ref{str-dis-2}. In Fig. \ref{str-dis-1} the result of $h_1$ on the 2-D triangular mesh is compared with the reference solution presented in \cite{kurganov2009-smallpertur}. Good agreement has been obtained. The water levels of the first layer together with the water surface and discharge are plotted in Fig. \ref{str-dis-2}.

\begin{figure}[!htb]
\centering
\includegraphics[width=0.40\textwidth]{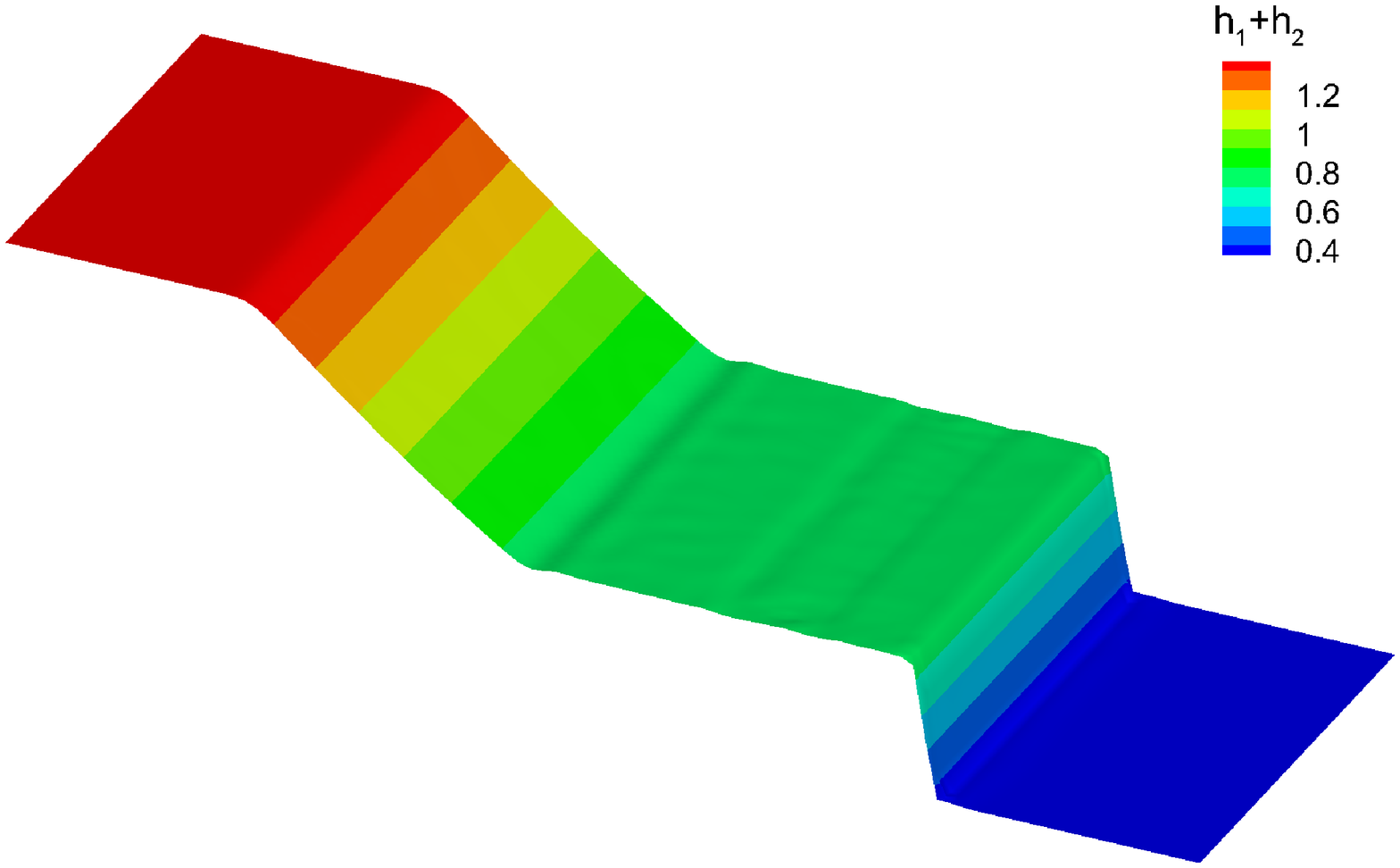}
\includegraphics[width=0.40\textwidth]{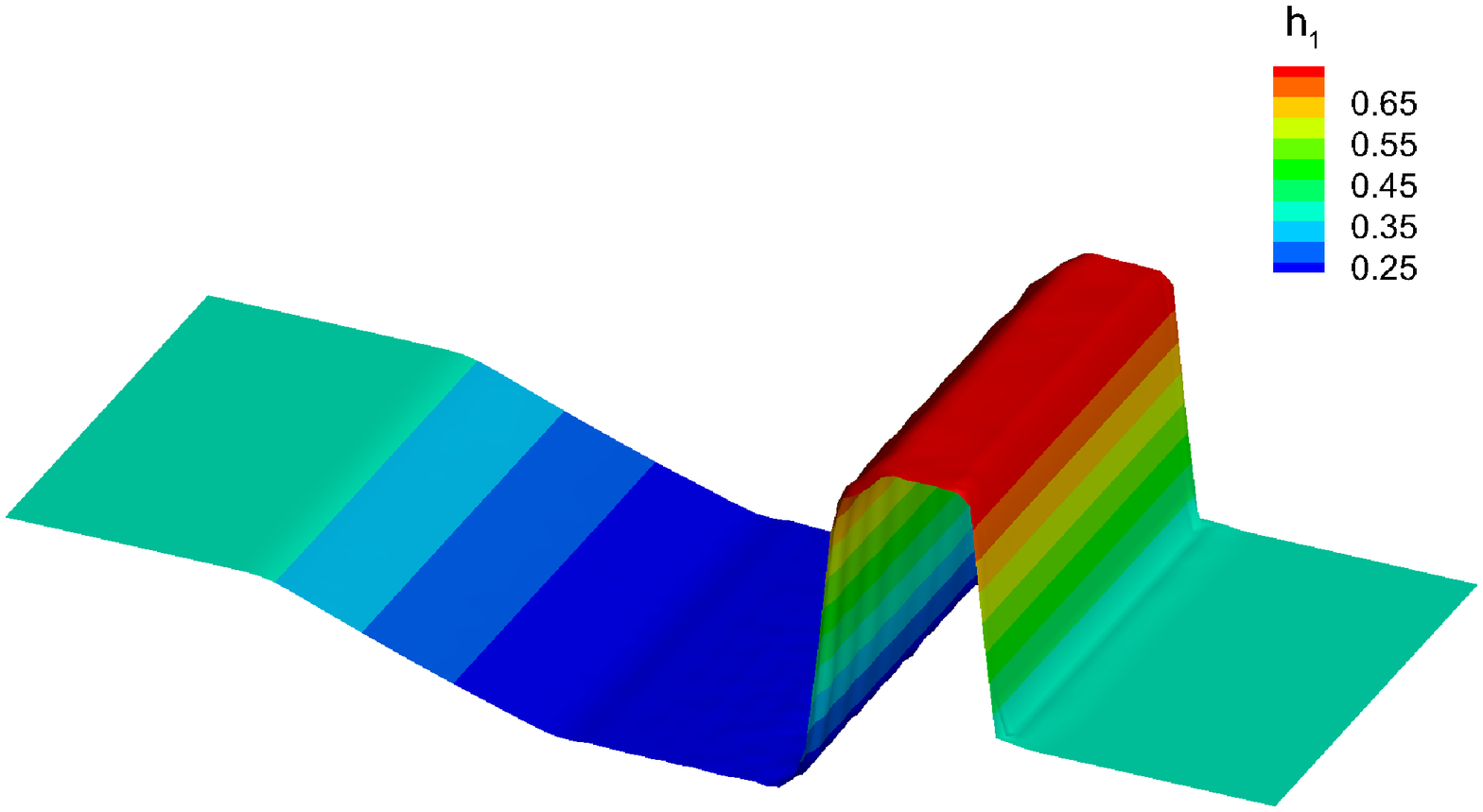}
\caption{\label{1d-dambreak-1-1} Dam-break flow at $\chi=1$: the 3-D contours of $h_1+h_2$ (left) and $h_1$ (right) at $t=0.08$. The cell sizes of the coarse and fine meshes are $\Delta X_{CM}=1/100$ and $\Delta X_{FM}=1/400$.}
\end{figure}

\begin{figure}[!htb]
\centering
\includegraphics[width=0.45\textwidth]{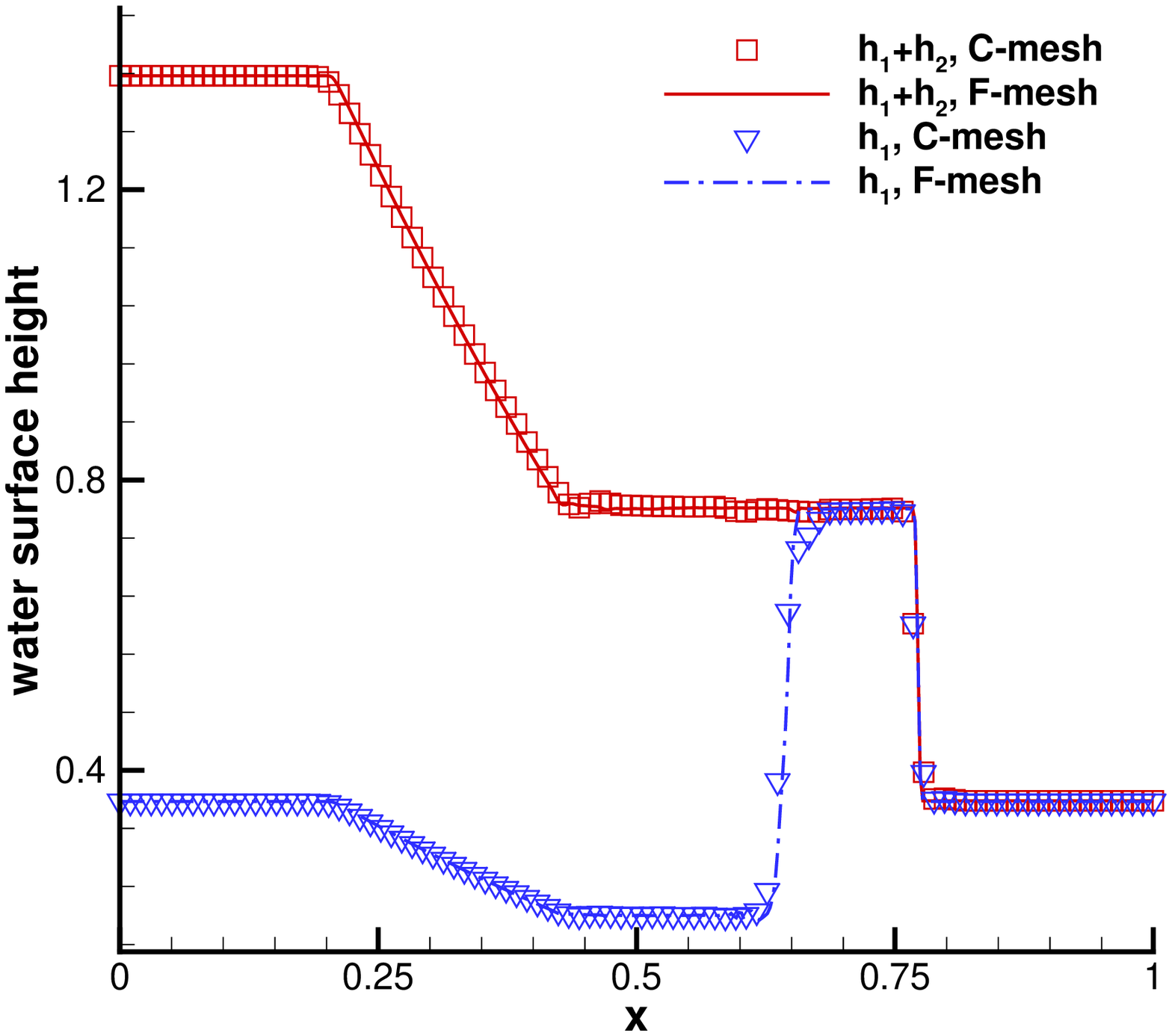}
\includegraphics[width=0.45\textwidth]{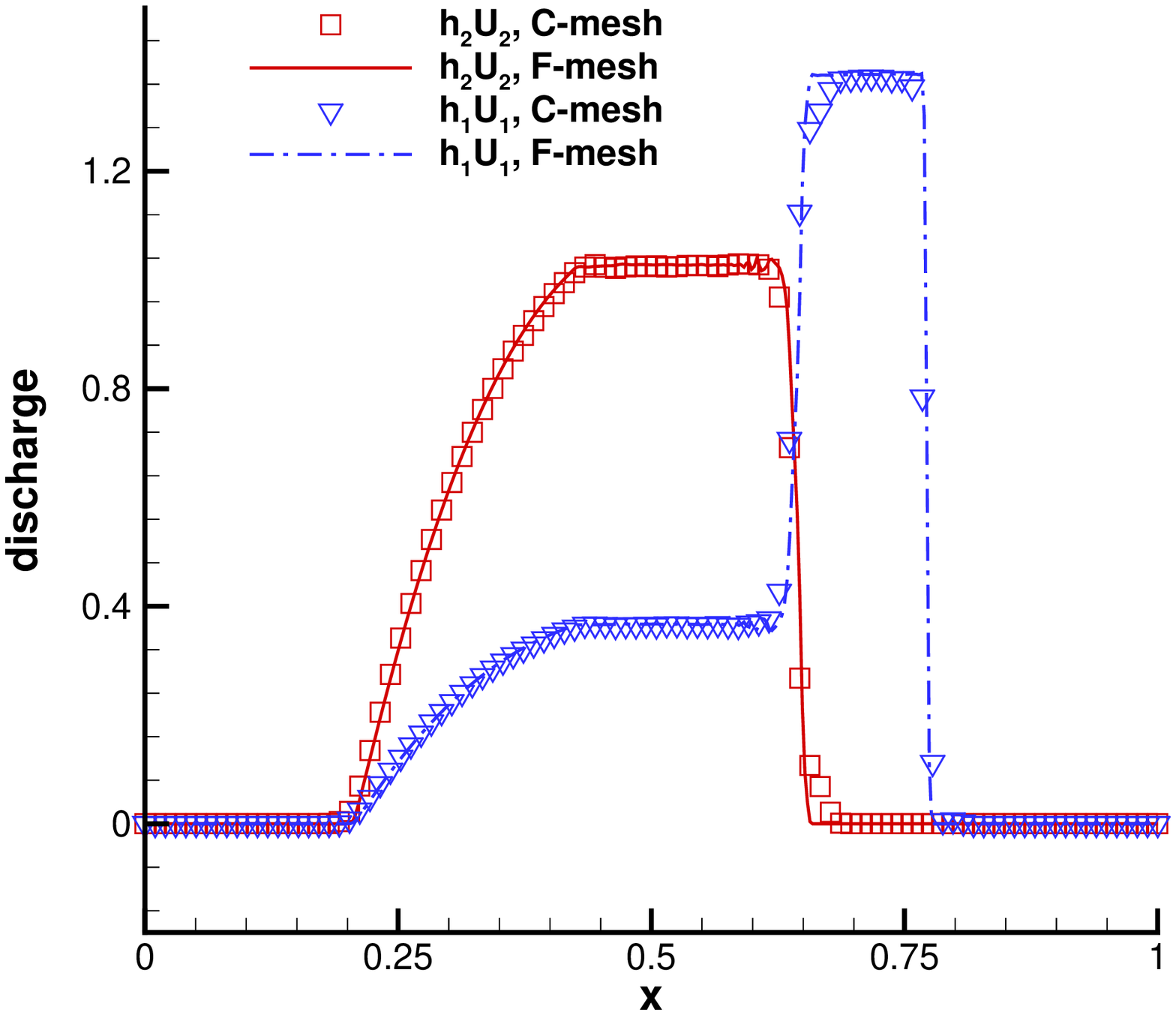}
\caption{\label{1d-dambreak-1-2} Dam-break flow at  $\chi=1$: the water levels (left)and discharge distributions (tight) along the horizontal centerline at $t=0.08$. The cell sizes of the coarse and fine meshes are $\Delta X_{CM}=1/100$ and $\Delta X_{FM}=1/400$.}
\end{figure}

\subsection{Dam-break problems at different density ratios}
The two-layer dam-break flows are used to validate the compact GKS.
The initial state is given as
\begin{equation*}
(h_1,h_2) = \begin{cases}
(0.357,1),  0\leq x<0.5,\\
(0.357,0), 0.5\leq x\leq1.
\end{cases}
\end{equation*}
The velocity is set as $(U_1,V_1)=(U_2,V_2)=(0,0)$ in the whole domain, and the computational domain is $[0,1]\times[0,0.5]$.
The gravitational acceleration is $G=9.81$. Dam-break flows at two density ratios are studied. In the computation, a coarse triangular mesh with $\Delta X=1/100$ and a fine triangular mesh with $\Delta X=1/400$ are used.

The first case is the dam-break flow at same density of the two layers, i.e., the density ratio with $\chi=1$.
The 3-D water level distributions of $h_1+h_2$ and $h_1$ at $t=0.08$ obtained by the compact GKS on the coarse mesh are shown in Fig. \ref{1d-dambreak-1-1}.
The water levels and discharge distributions along the horizontal centerline are given in Fig. \ref{1d-dambreak-1-2}.
The results on the coarse mesh are consistent with those on the fine mesh, and the water levels obtained by the current compact scheme are consistent with those in \cite{TLSWE-3forms}.

\begin{figure}[!htb]
\centering
\includegraphics[width=0.40\textwidth]{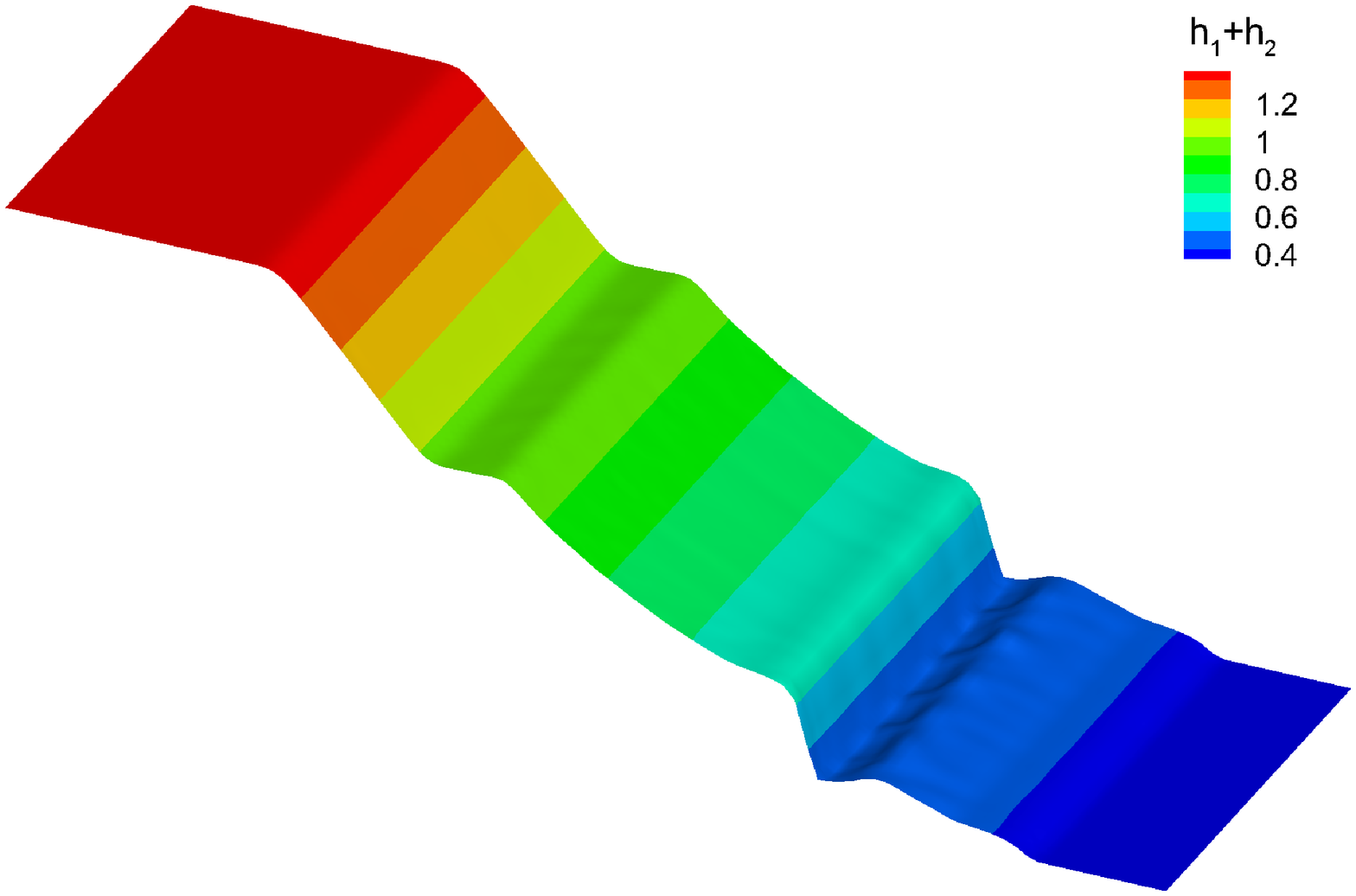}
\includegraphics[width=0.40\textwidth]{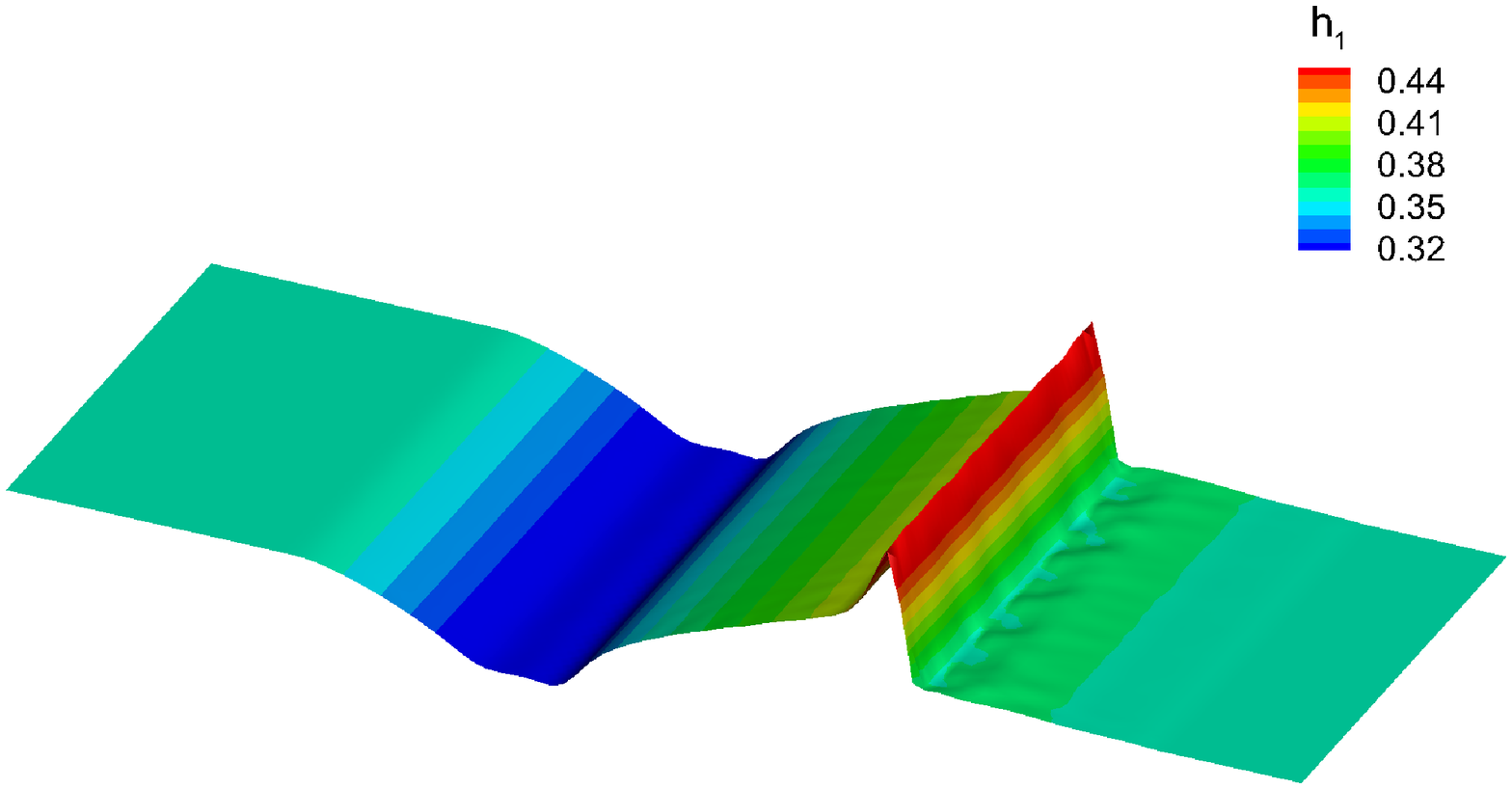}
\caption{\label{1d-dambreak-2-1} Dam-break flow of a light fluid over a dense fluid with $\chi=0.2$: the 3-D height distributions of $h_1+h_2$ (left) and $h_1$ (right) at $t=0.08$. The cell sizes of the coarse and fine meshes are $\Delta X_{CM}=1/100$ and $\Delta X_{FM}=1/400$. }
\end{figure}

\begin{figure}[!htb]
\centering
\includegraphics[width=0.45\textwidth]{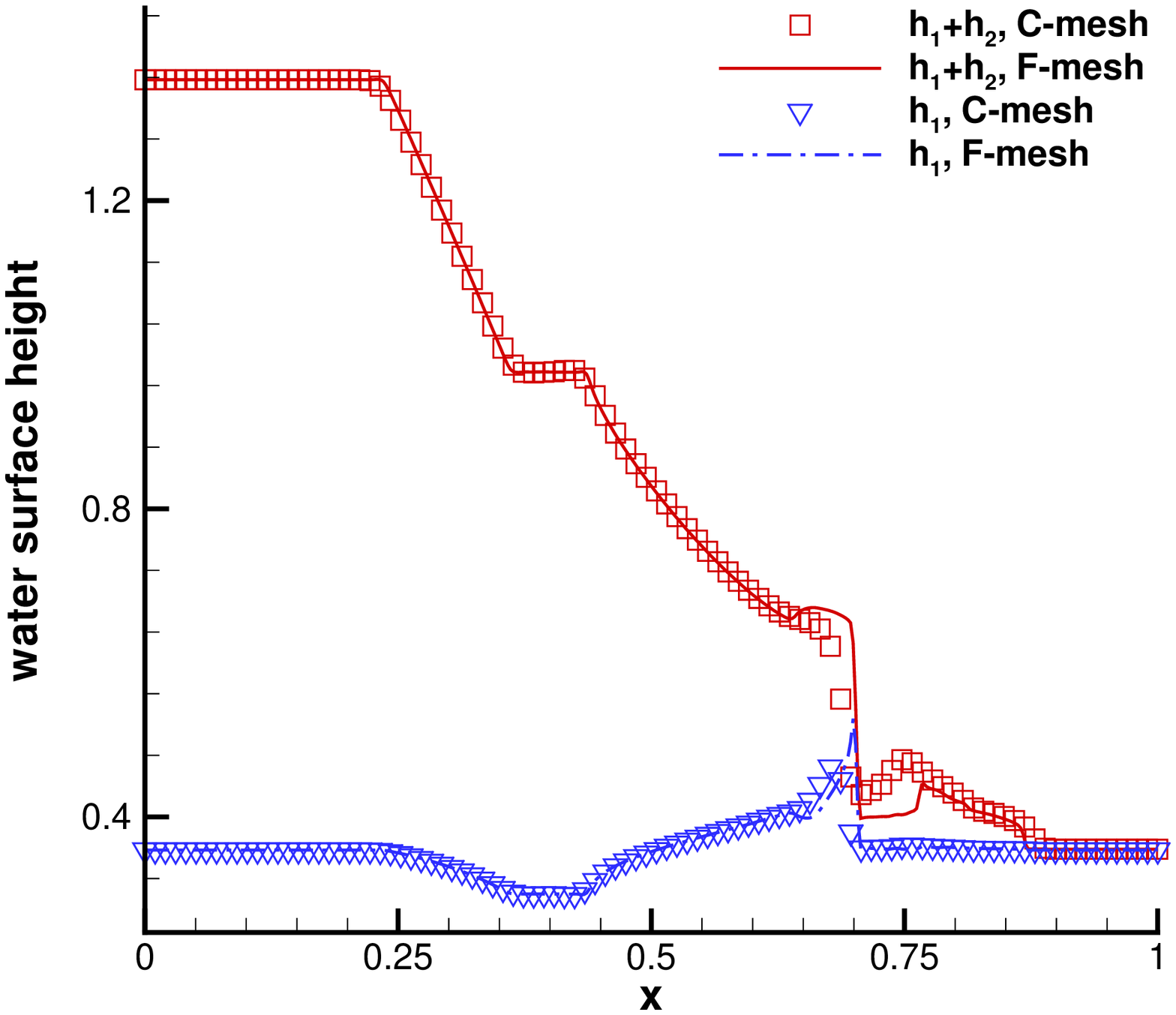}
\includegraphics[width=0.45\textwidth]{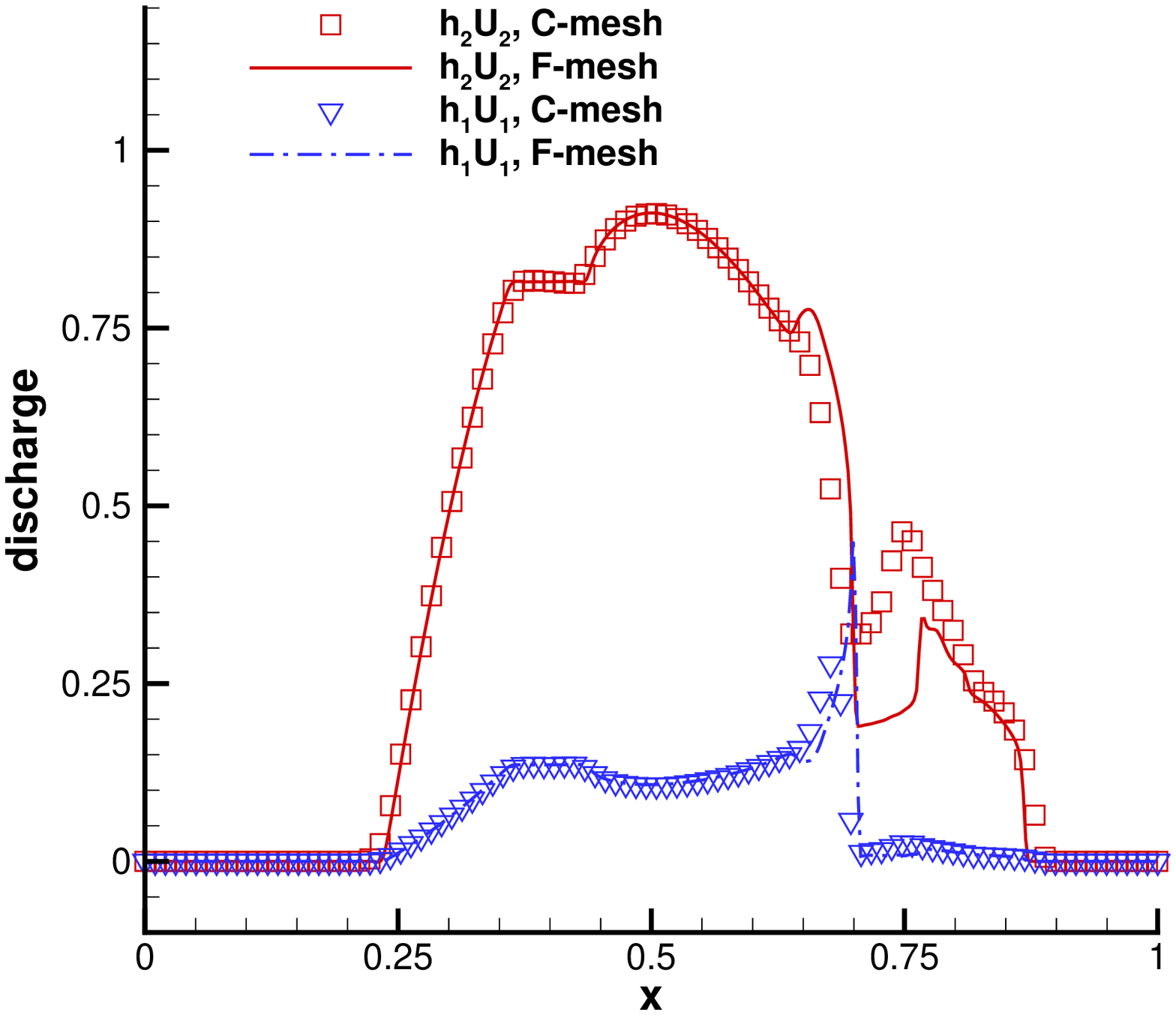}
\caption{\label{1d-dambreak-2-2} Dam-break flow of a light fluid over a dense fluid with $\chi=0.2$: the height (left) and discharge (right) along the horizontal centerline at $t=0.08$. The cell sizes of the coarse and fine meshes are $\Delta X_{CM}=1/100$ and $\Delta X_{FM}=1/400$. }
\end{figure}

The second case is the dam-break flow of a light fluid over a dense one. The density ratio is $\chi=0.2$.
The 3-D water level distributions of $h_1+h_2$ and $h_1$ at $t=0.08$ obtained by the compact GKS are shown in Fig.\ref{1d-dambreak-2-1}.
The 1-D water levels and discharge distributions along the horizontal centerline are given in Fig.\ref{1d-dambreak-2-2}.
Due to the complexity of the solution, the fine mesh result has a better spatial resolution and gives the solution close to the reference ones in \cite{TLSWE-3forms}.

\begin{figure}[!htb]
\centering
\includegraphics[width=0.495\textwidth]{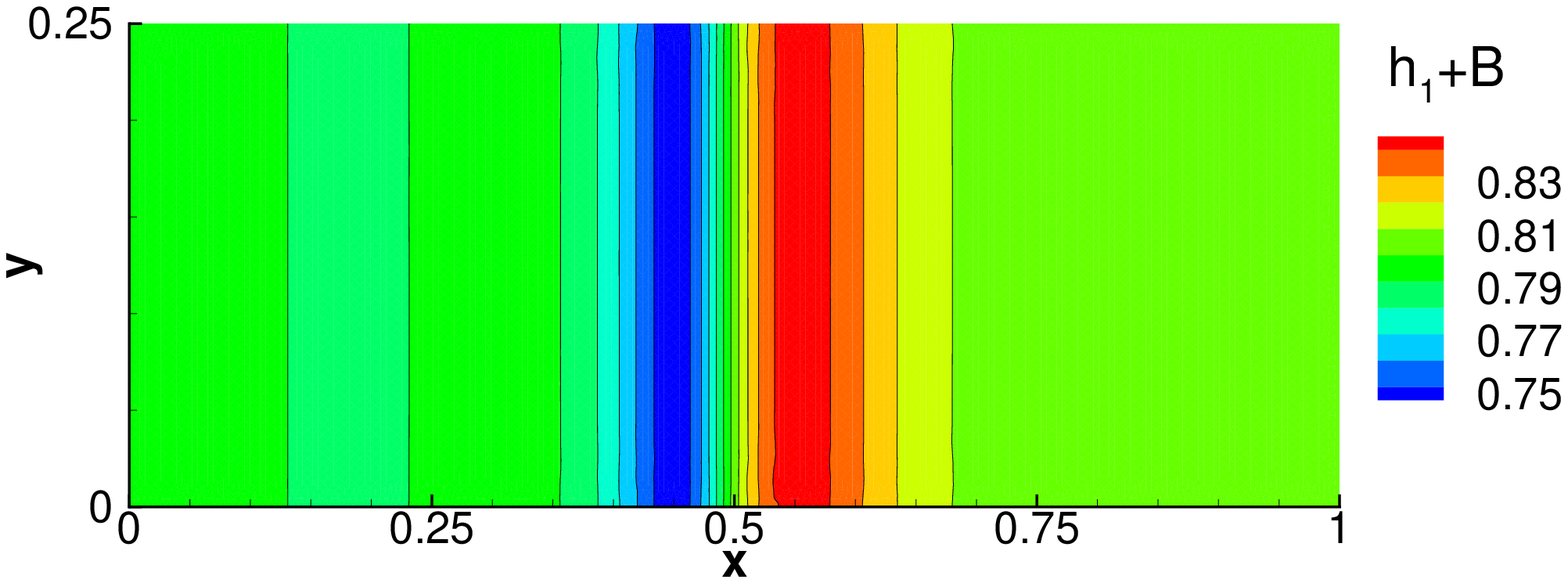}
\includegraphics[width=0.495\textwidth]{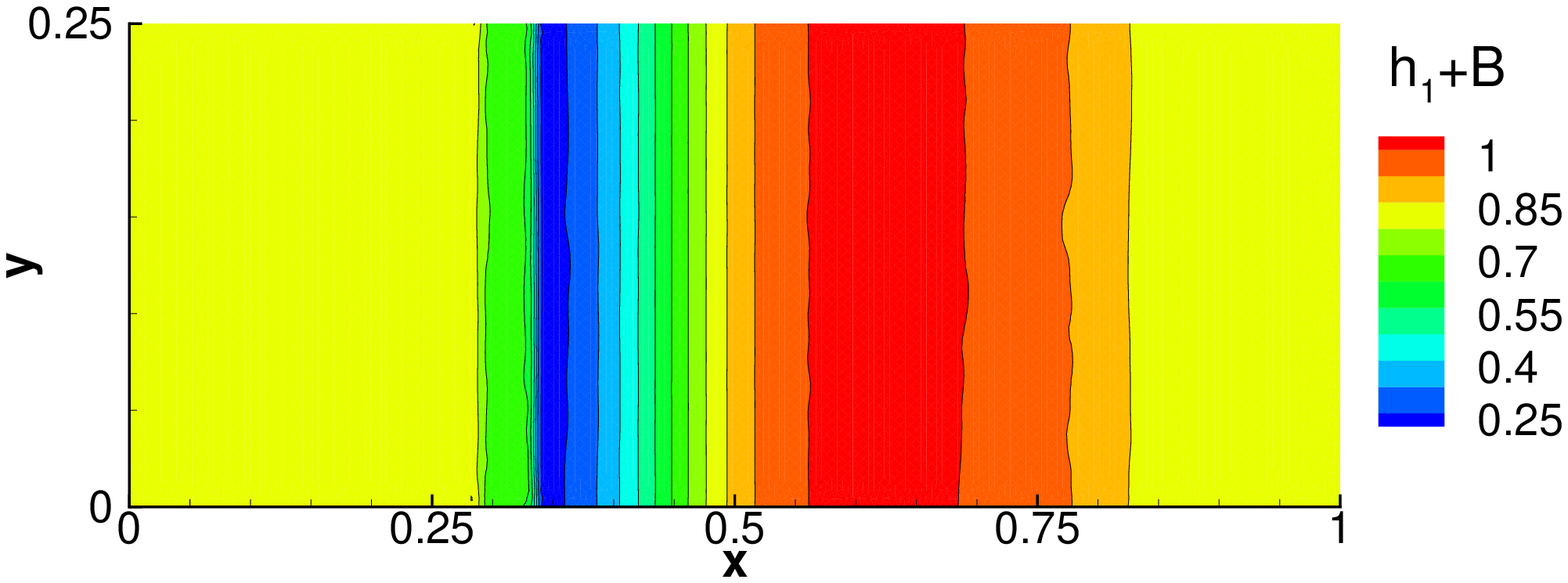}\\
\includegraphics[width=0.495\textwidth]{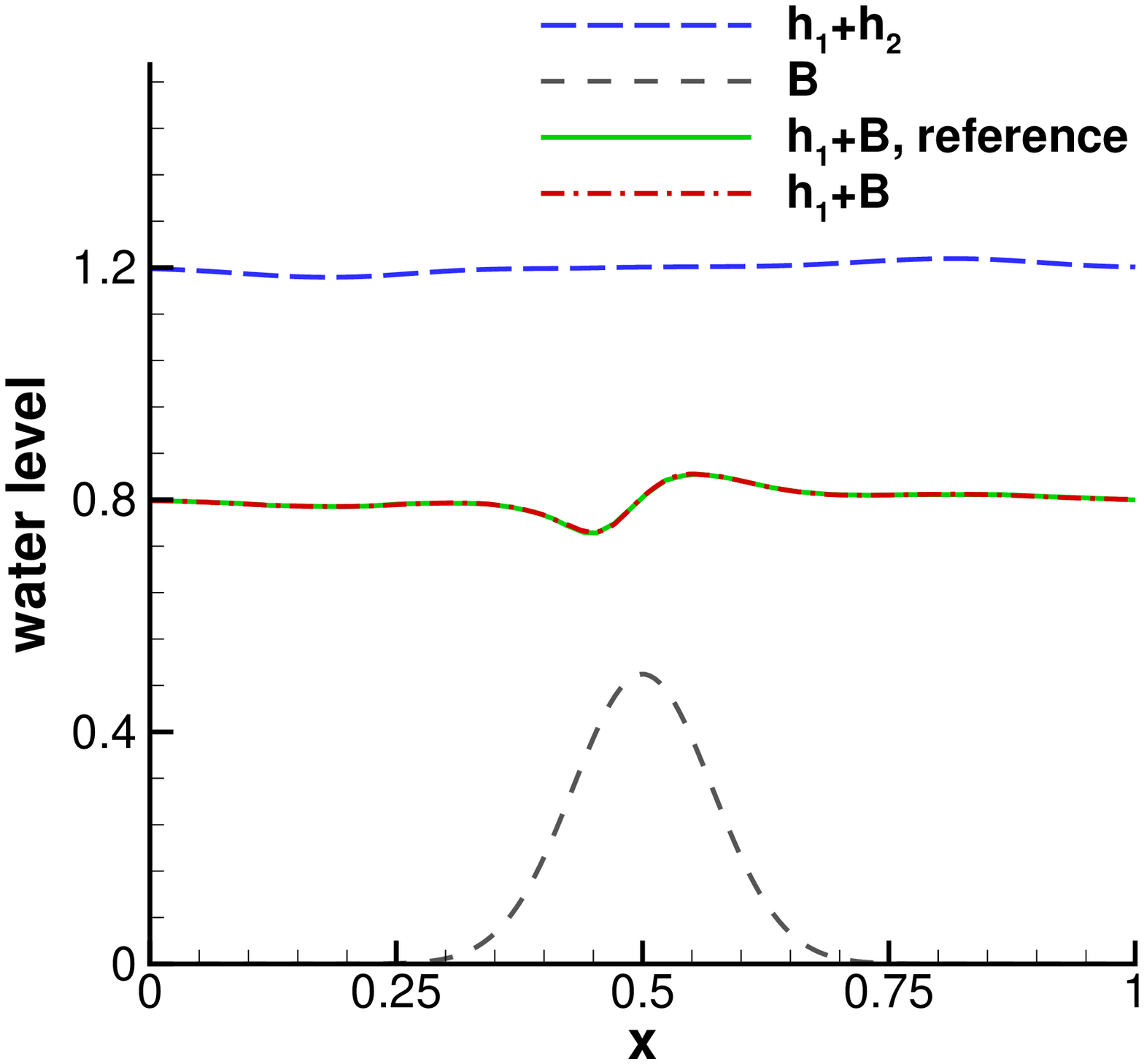}
\includegraphics[width=0.495\textwidth]{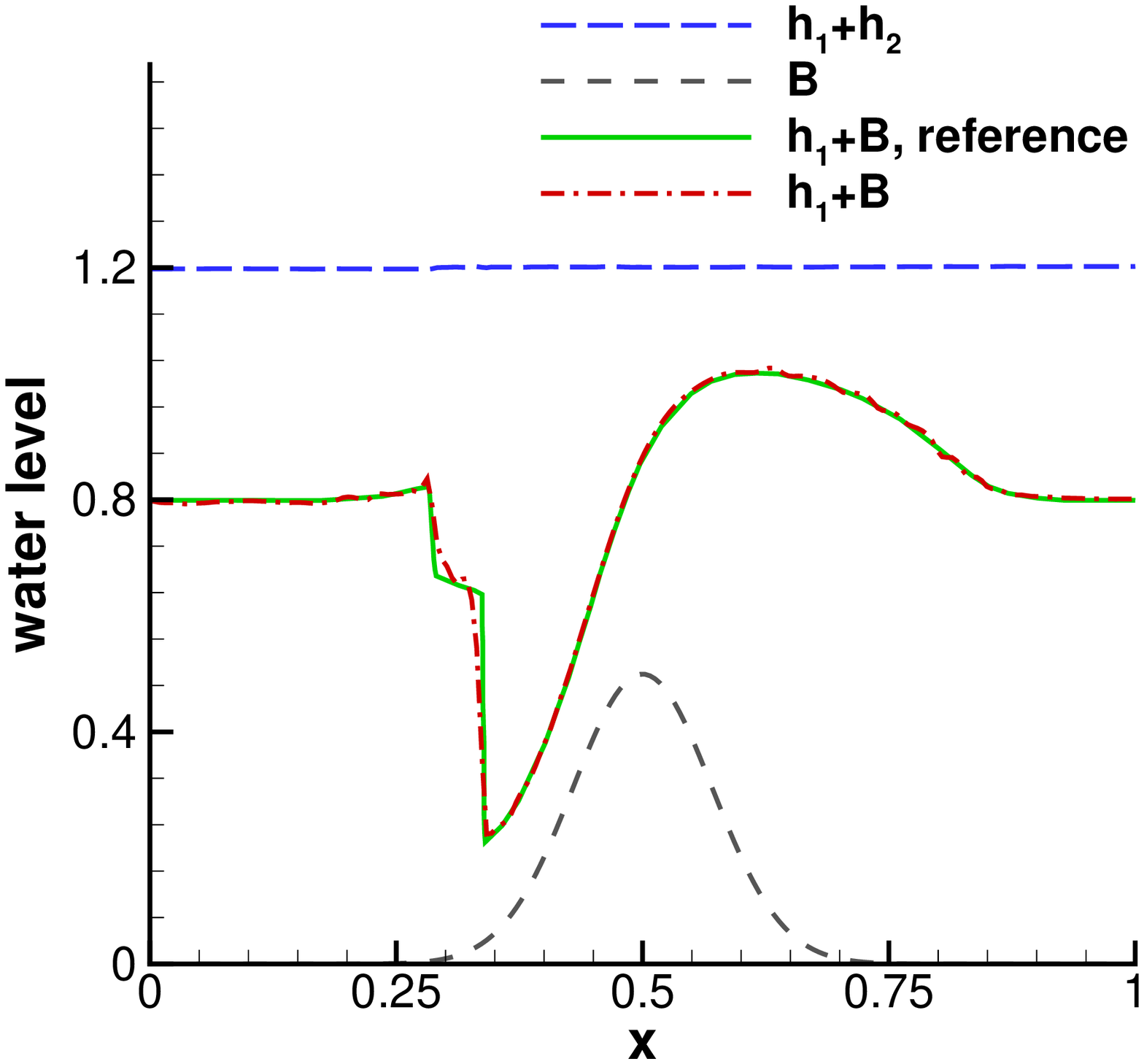}
\caption{\label{1d-nonflat} Channel flow with non-flat bottom:  the 2-D contours of the water levels (up) and the water level distributions along the horizontal centerline (down) at $t=0.1$ (left) and $t=1.0$ (right).}
\end{figure}

\subsection{Channel flow with non-flat bottom}

This case is about the two-layer flow through a channel with non-flat bottom topography.
The bottom topography is defined by
\begin{equation*}
B(x,y)=0.5e^{-100(x-0.5)^2}.
\end{equation*}
The initial condition is given as
\begin{equation*}
\begin{split}
h_1&=0.8-B(x,y),~~h_2=0.4, \\
U_1&=-0.2,~~~~~~~~~~~~U_2=0.15.
\end{split}
\end{equation*}
The channel covers a domain $[0,1]\times[0,0.25]$.
The reflecting boundary condition is applied at the channel walls.
The free boundary condition is used on the left and
right boundaries.
The triangular mesh with a cell size of $\Delta X=1/200$ is used in the computation.

Fig. \ref{1d-nonflat} shows the results of water levels at $t=0.1$ and $t=1.0$, respectively.
Due to the non-flat bottom topography, the interface between two layer fluids evolves from an initial smooth interface to a discontinuous one.
The reference solution comes from solving the 1-D TLSWE on a uniform mesh with $1000$ cells in \cite{2018nonflat-channel}.
At the early time, a smooth interface evolves, such as the left figures in Fig. \ref{1d-nonflat}, and the solution has good agreement with the reference solution. At a later time, a discontinuous interface emerges, such as the right figures in Fig. \ref{1d-nonflat},
and the position of the discontinuity obtained by the compact GKS has a good match with the reference solution.

\subsection{2-D interface propagation}
The 2-D circular  interface propagation Riemann problem is studied.
The initial condition of the test case is given by
\begin{equation*}
(h_1,h_2) = \begin{cases}
(1.8,0.2), ~~ (x-5)^2+(y-5)^2 <4.0,\\
(0.2,1.8), ~~ \mathrm{otherwise}.
\end{cases}
\end{equation*}
The initial velocity is $(U_1,V_1)=(U_2,V_2)=(0,0)$ in the computational domain $[0,10]\times[0,10]$.
The gravitational acceleration is $G=9.81$. The density ratio between layers is $\chi=0.98$.
The free boundary condition is adopted on all boundaries. The triangular mesh with a cell size of $\Delta X=1/10$ is used in the computation.

\begin{figure}[!htb]
\centering
\includegraphics[width=0.32\textwidth]{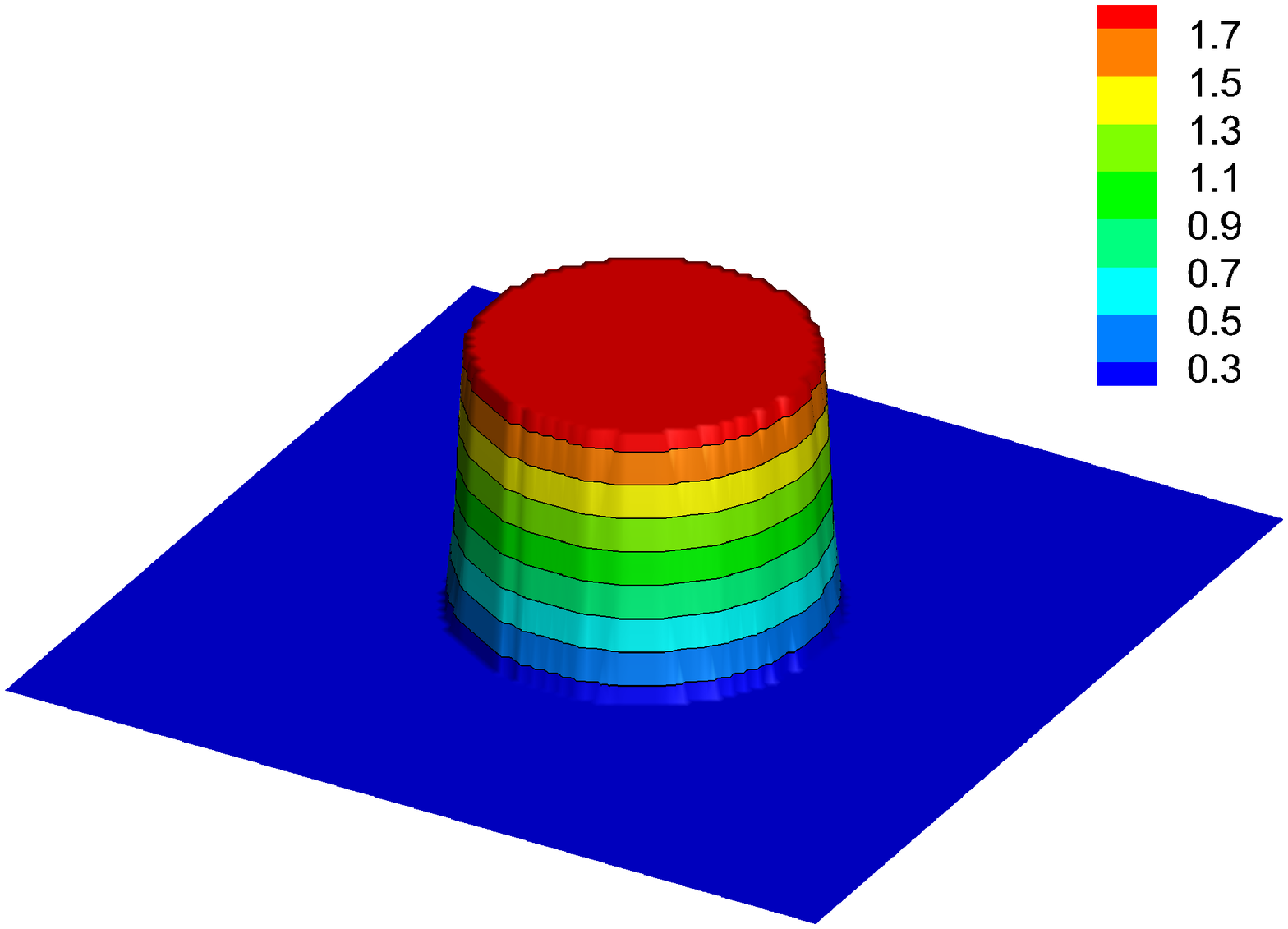}
\includegraphics[width=0.32\textwidth]{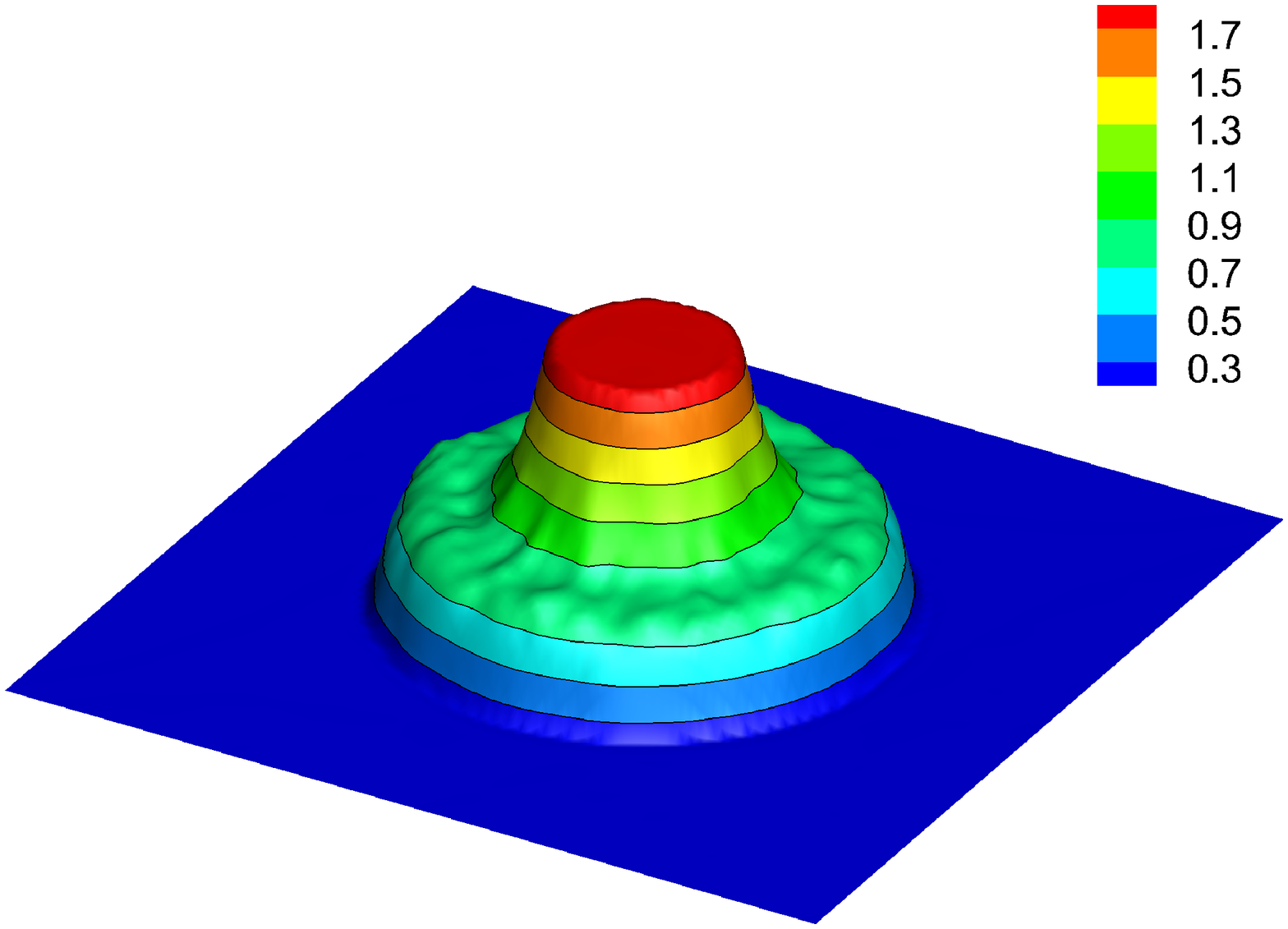}
\includegraphics[width=0.32\textwidth]{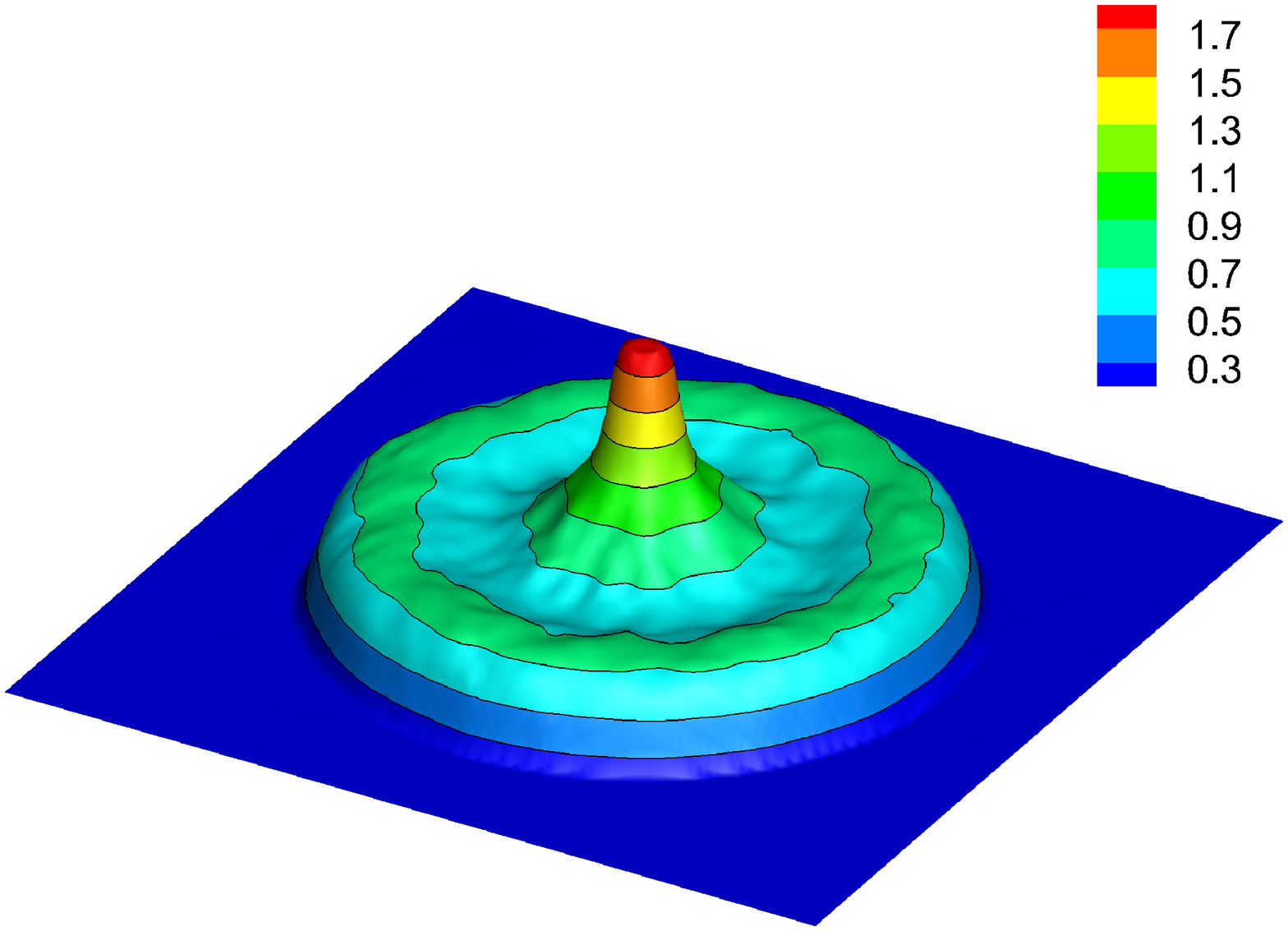}
\caption{\label{2d-dambreak-circular-1} 2D interface propagation: the 3-D water level distributions of $h_1$ at $t=0$, $t=2.0$ and $t=4.0$. The cell size of the triangular mesh is $\Delta X=1/10$. }
\end{figure}

The 3-D water level distributions of $h_1$ and its distributions along the horizontal centerline at $t=0$, $t=2.0$ and $t=4.0$ are presented in Fig. \ref{2d-dambreak-circular-1} and Fig. \ref{2d-dambreak-circular-1}, respectively.
The results show the circular propagation of the water column.

\begin{figure}[!htb]
\centering
\includegraphics[width=0.32\textwidth]{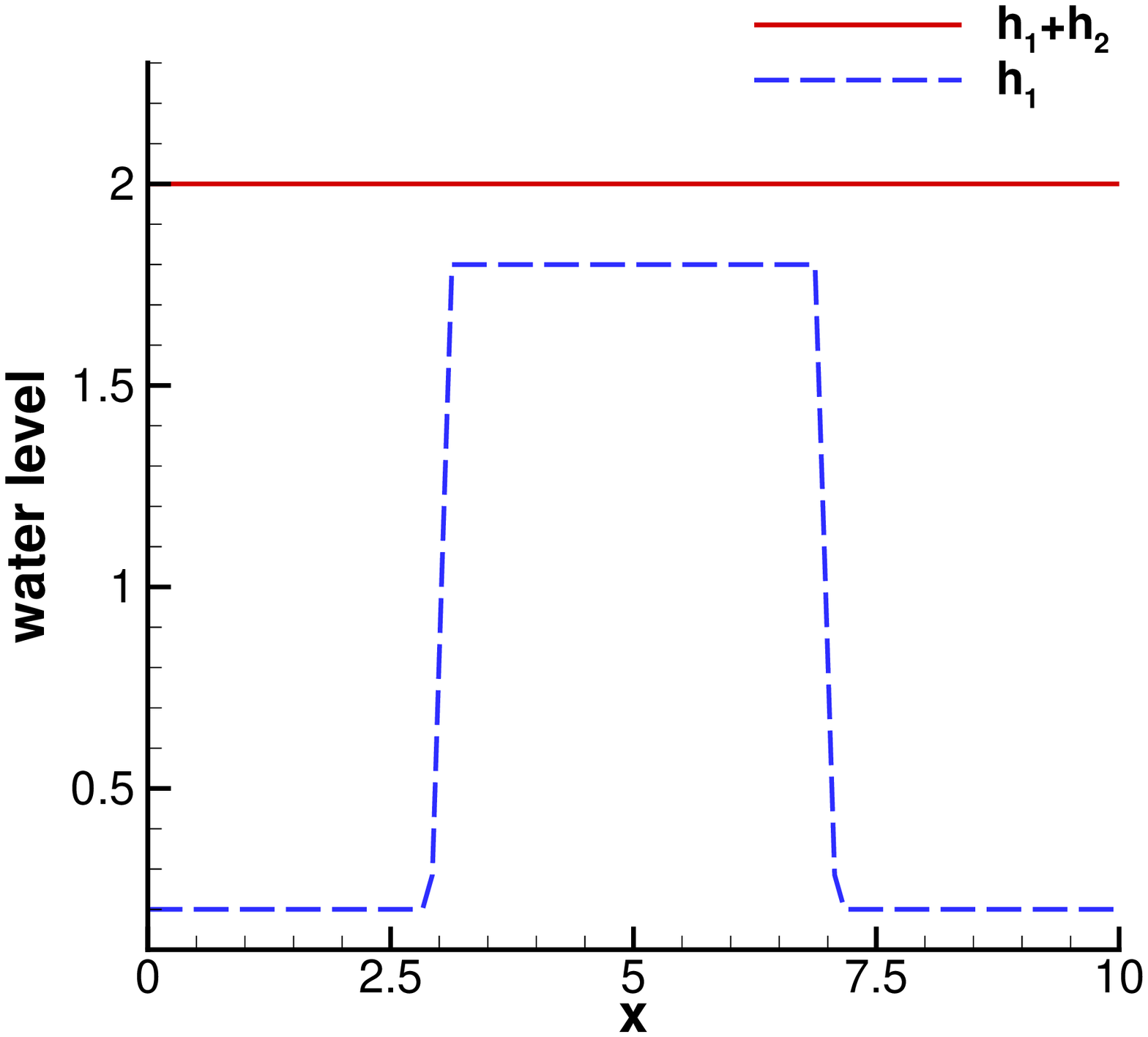}
\includegraphics[width=0.32\textwidth]{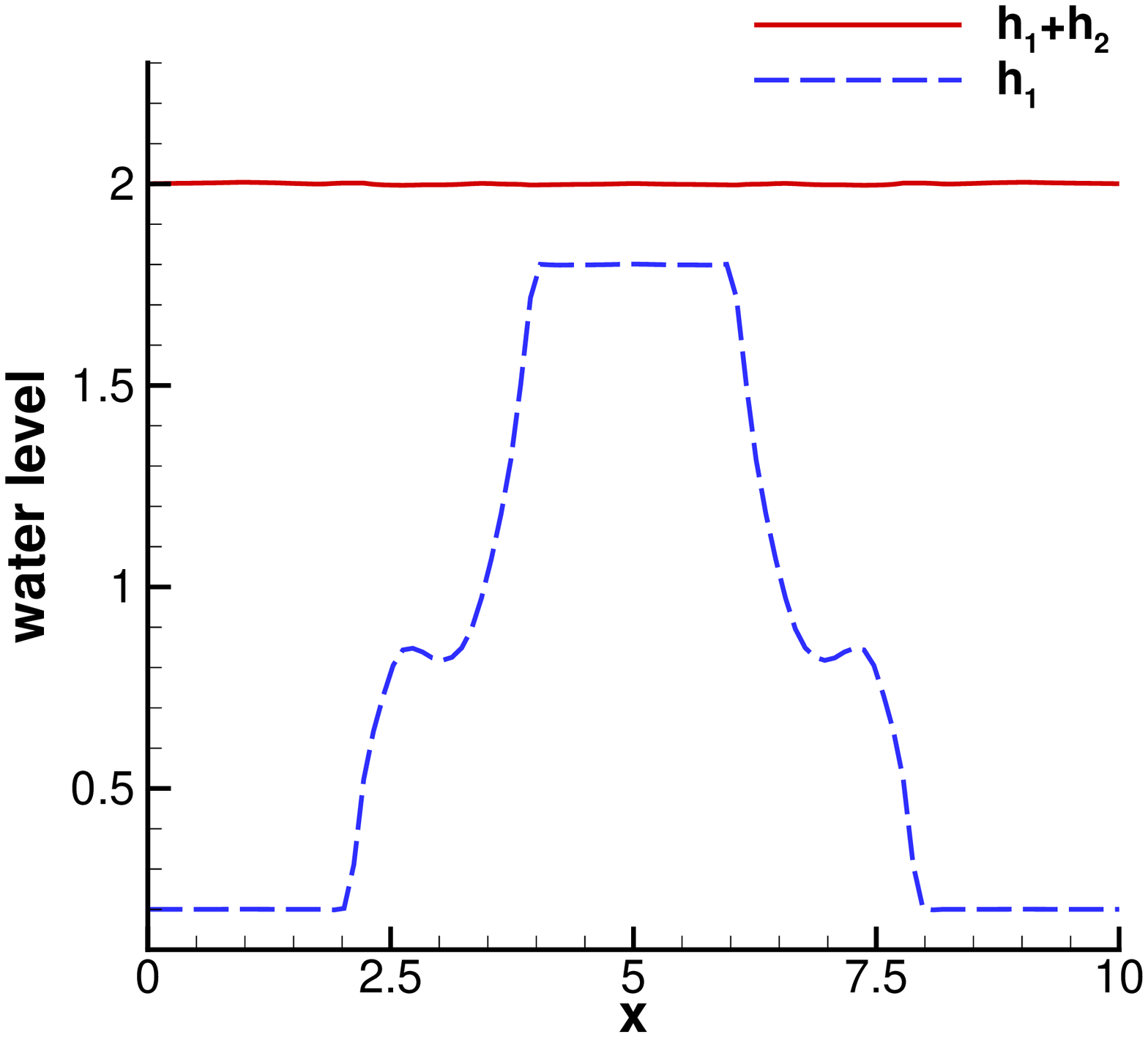}
\includegraphics[width=0.32\textwidth]{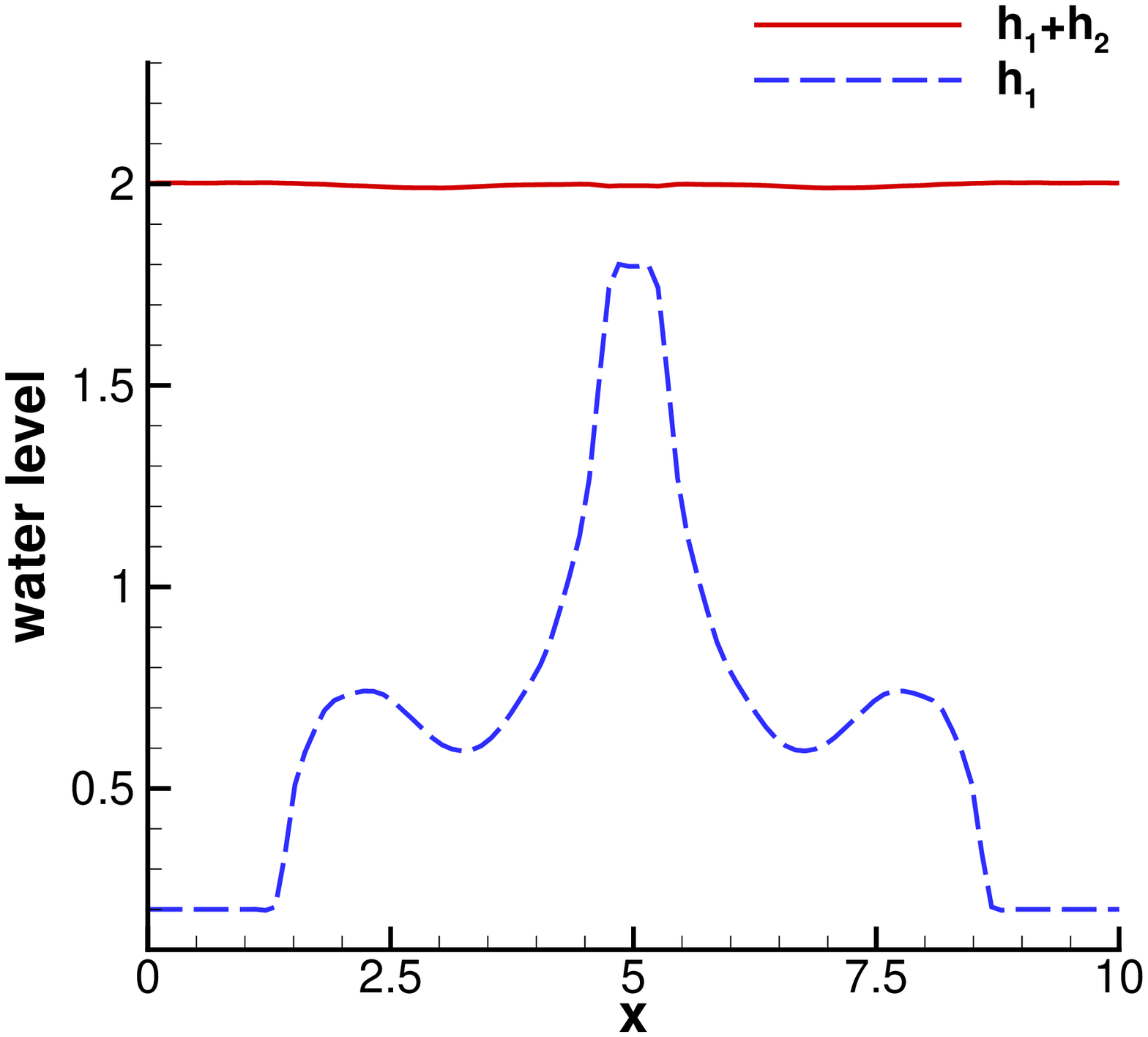}
\caption{\label{2d-dambreak-circular-2} 2D interface propagation: the water level distribution of $h_1$ along the horizontal centerline at $t=0$, $t=2.0$ and $t=4.0$. The cell size of the triangular mesh is $\Delta X=1/10$. }
\end{figure}

\begin{figure}[!htb]
	\centering
    \includegraphics[width=0.495\textwidth]{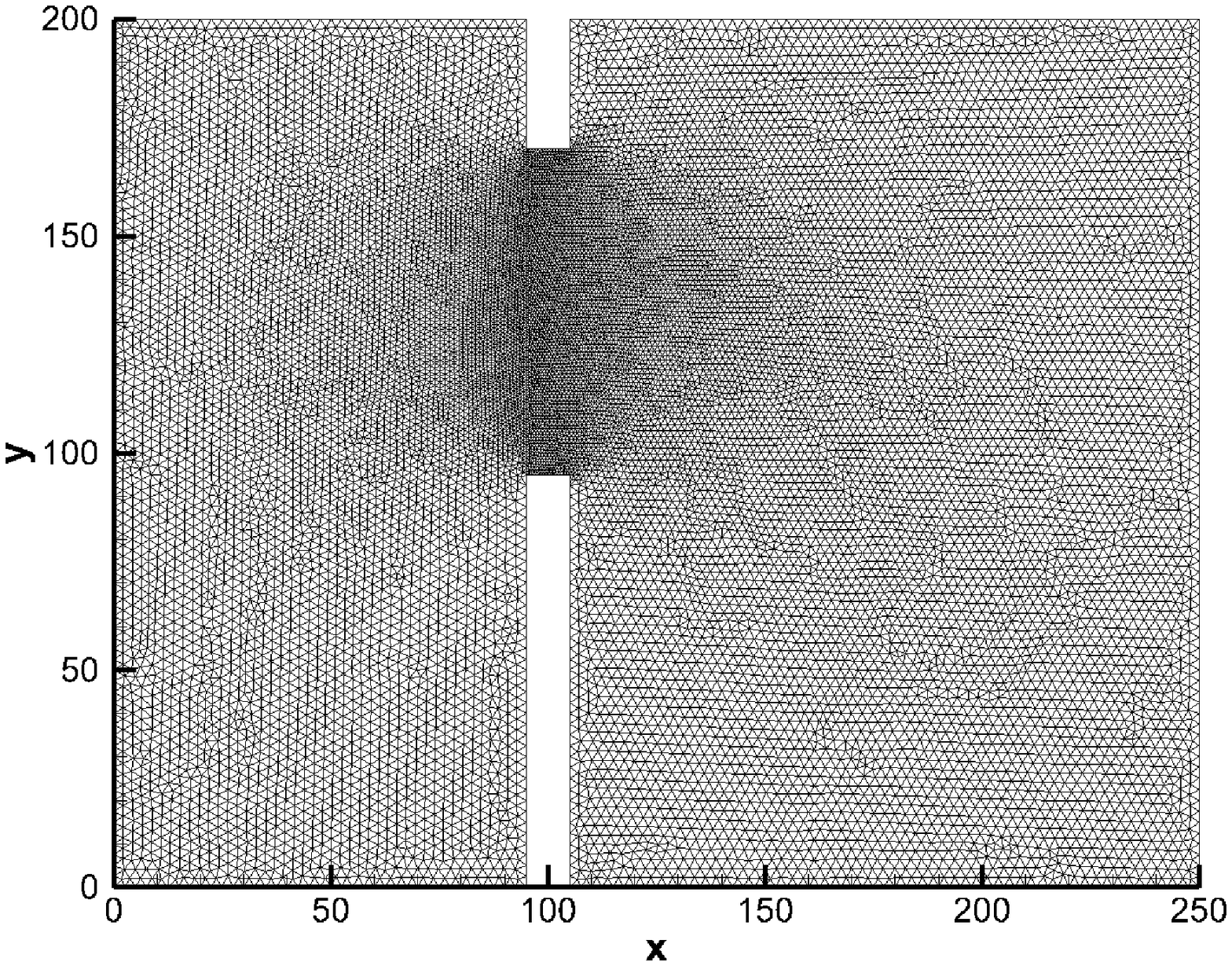}
    \includegraphics[width=0.495\textwidth]{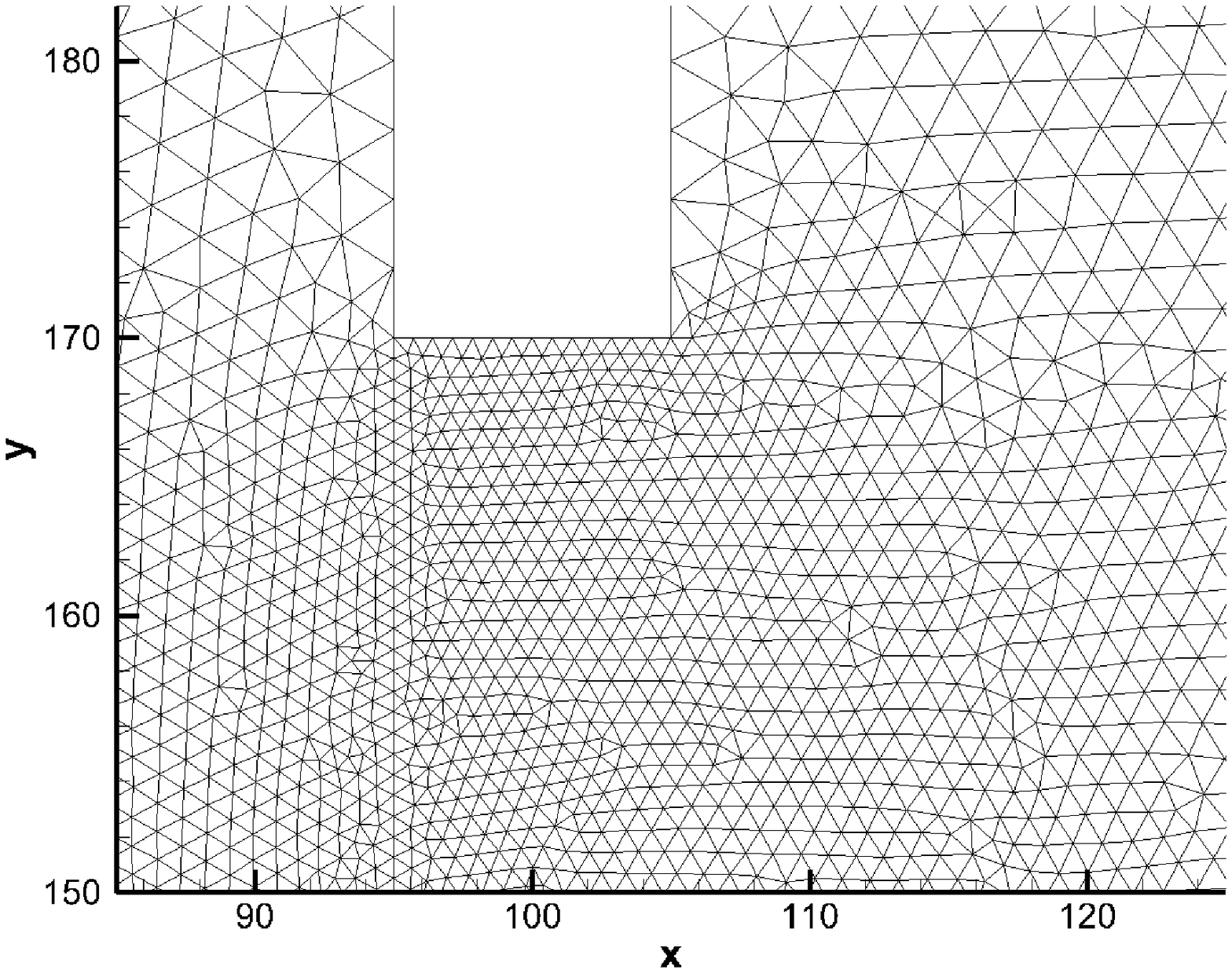}
	\caption{\label{2d-dambreak-mesh} 2-D dam-break in an irregular domain: the left figure is the computational domain and mesh, and the right figure is the enlarged view of the mesh around the dam breach. The mesh size far away from the dam is $\Delta X=2.5$, and the mesh size
is refined by $3.3$ times in the region close to dam breach. }
\end{figure}

\begin{figure}[!htb]
\centering
\includegraphics[width=0.32\textwidth]{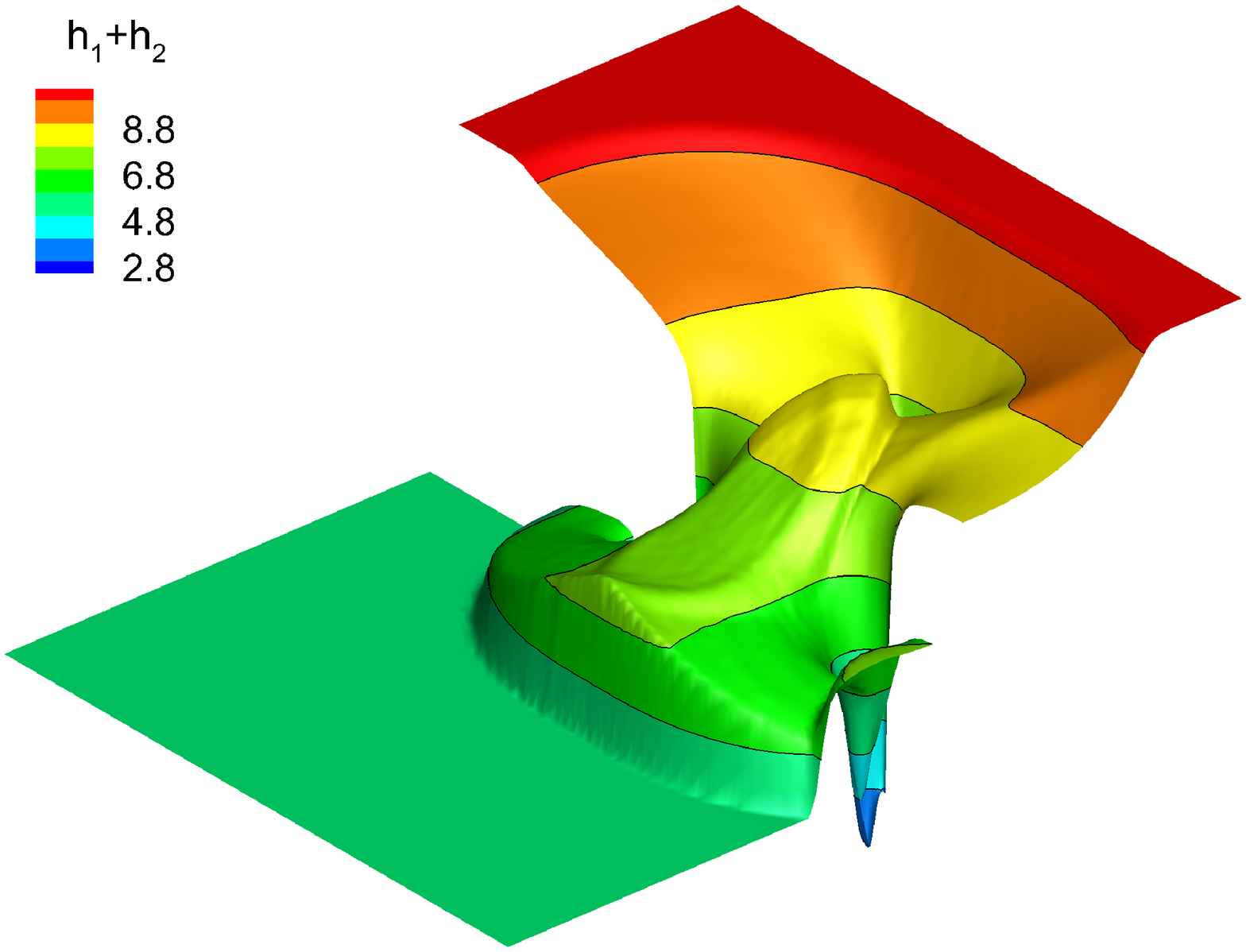}
\includegraphics[width=0.32\textwidth]{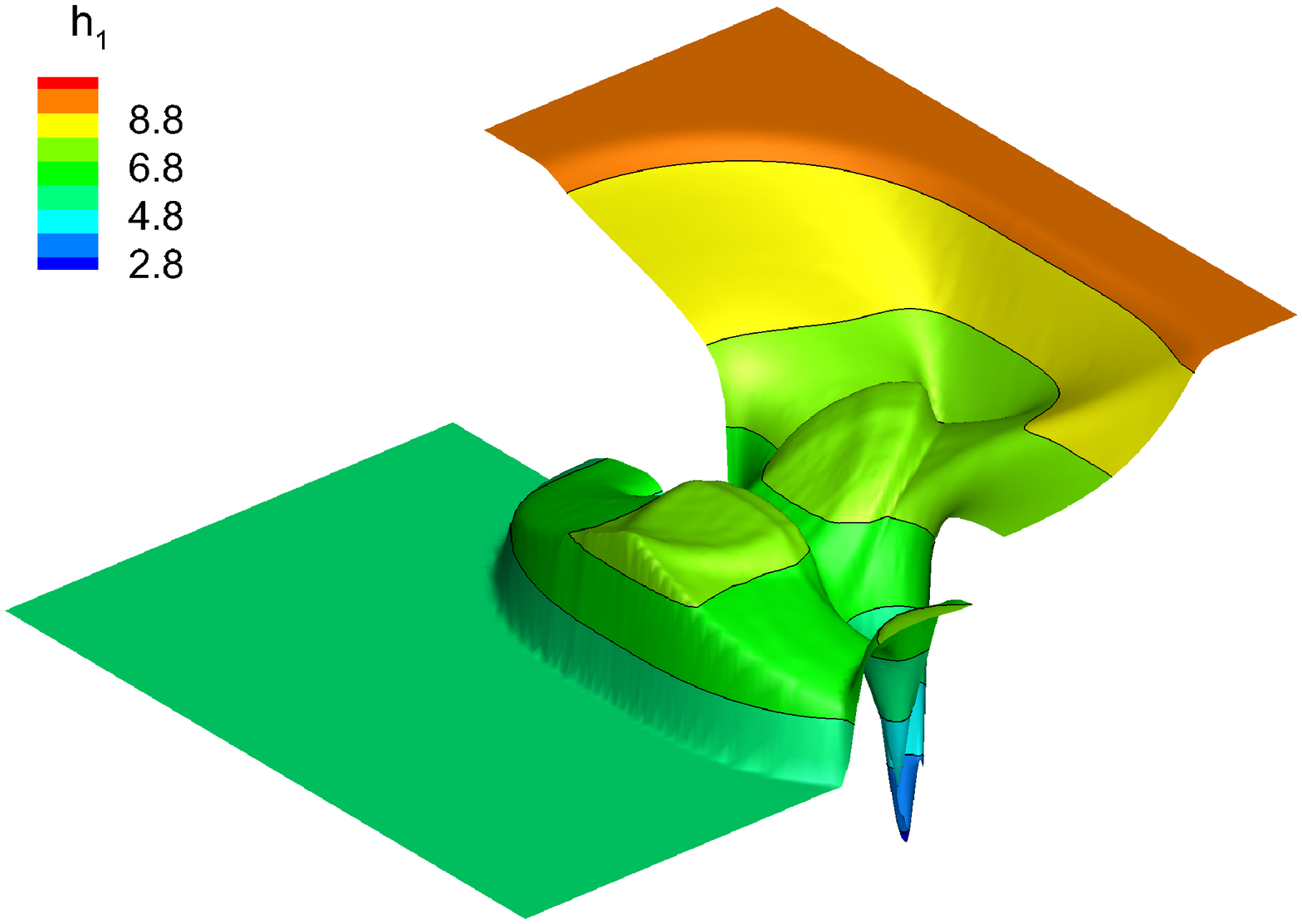}
\includegraphics[width=0.32\textwidth]{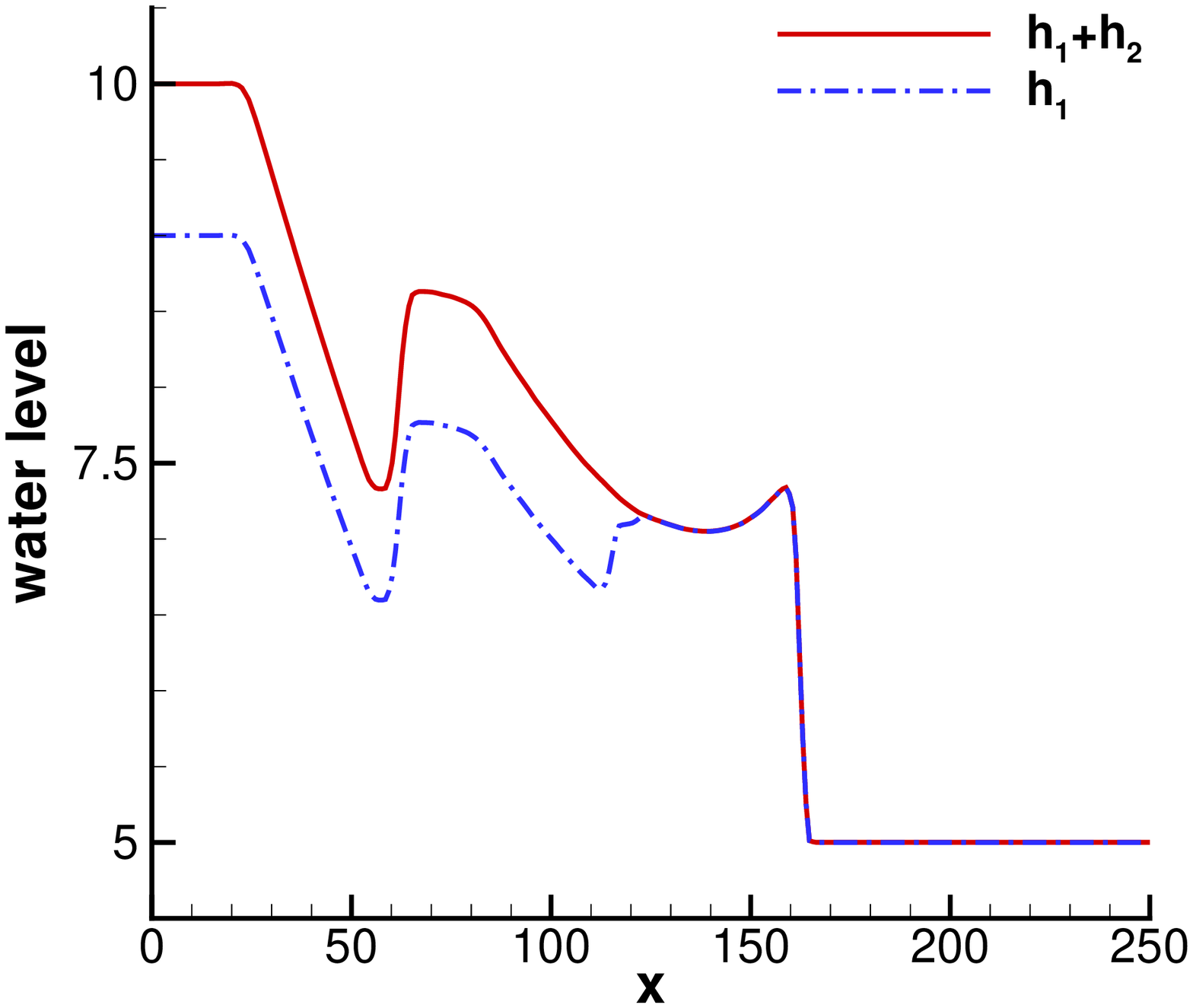}
\caption{\label{2d-dambreak-wet} 2-D dam-break in an irregular domain with a wet bed: the 3-D contours of water levels and the distributions along the horizontal centerline of the breach. The computational time is $t=7.2$.}
\end{figure}

\begin{figure}[!htb]
\centering
\includegraphics[width=0.32\textwidth]{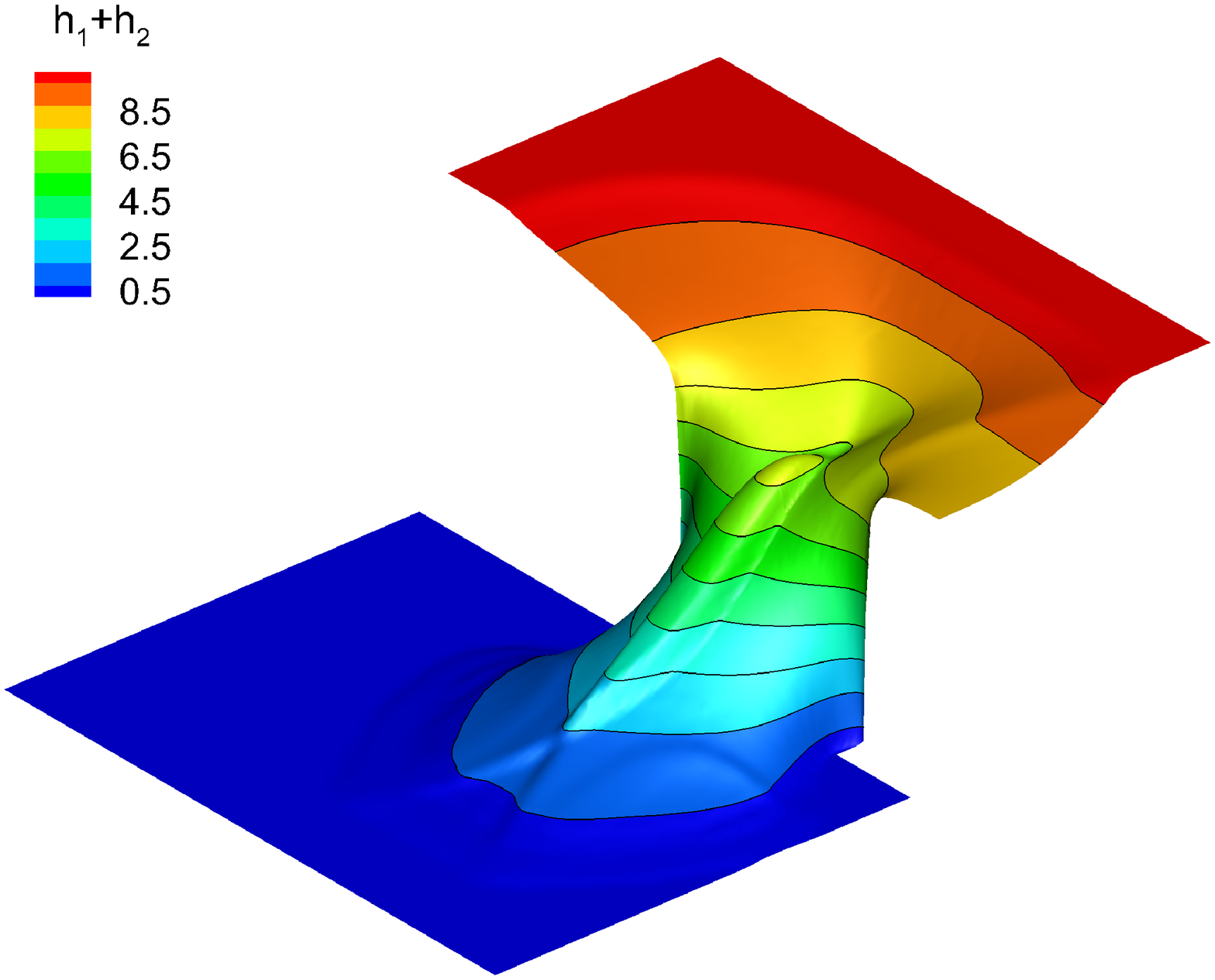}
\includegraphics[width=0.32\textwidth]{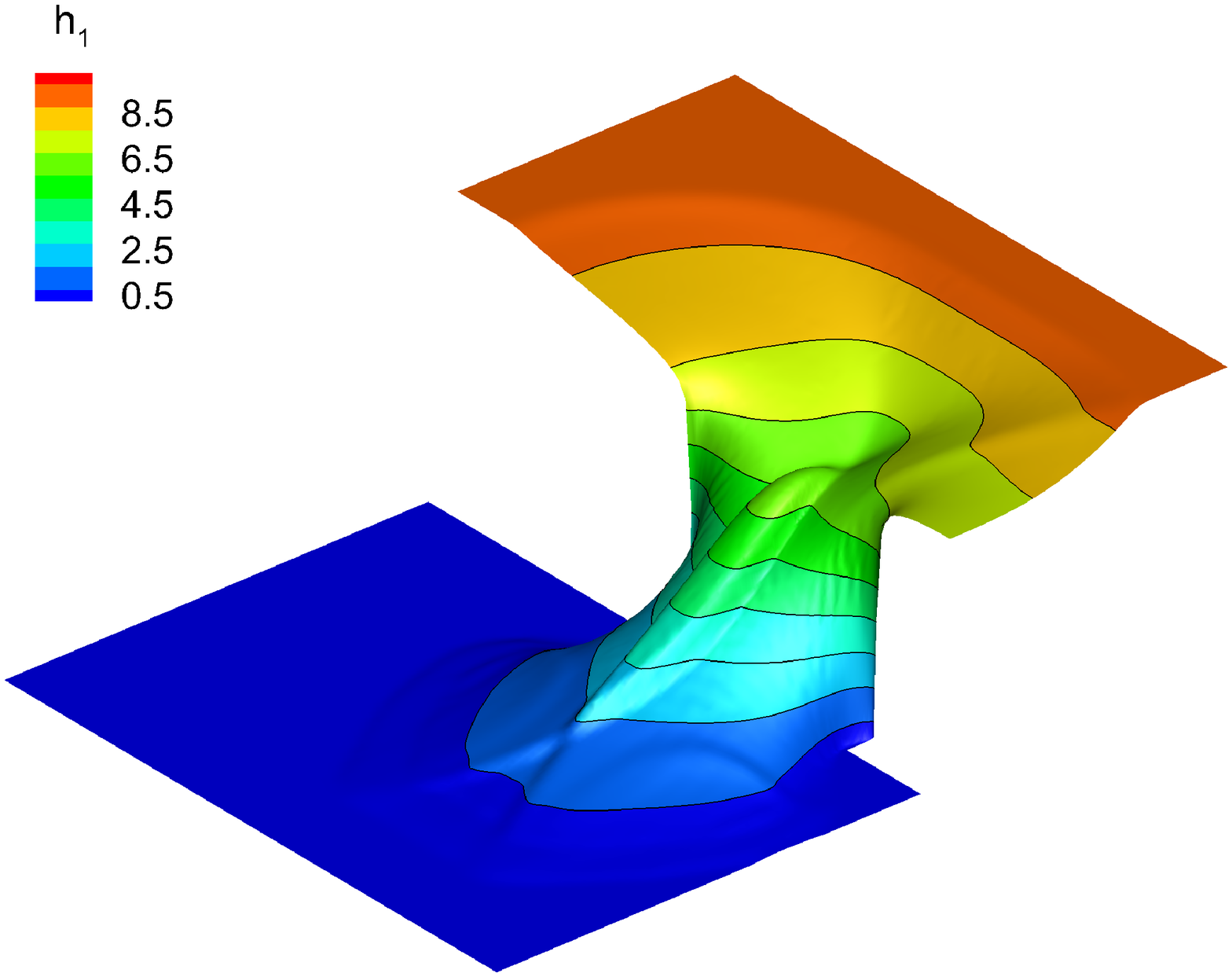}
\includegraphics[width=0.32\textwidth]{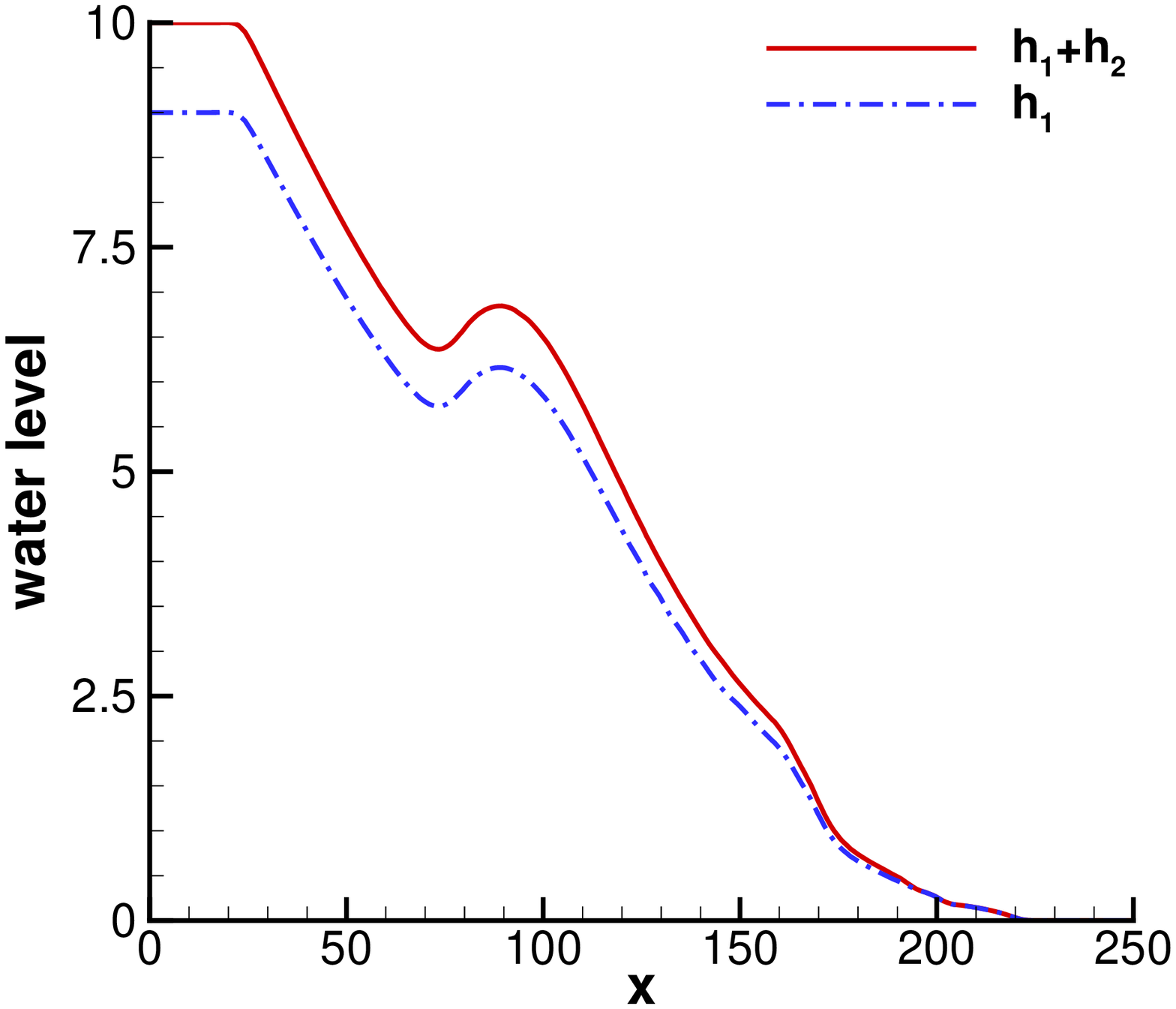}
\caption{\label{2d-dambreak-dry} 2-D dam-break in an irregular domain with a dry bed: the 3-D contours of water levels and the distributions along the horizontal centerline of the breach. The computational time is $t=7.2$.}
\end{figure}

\subsection{2-D dam-break in an irregular domain}
The 2-D dam-break problem in \cite{dambreak-1998,dambreak-2007} is used in the current study to validate the compact GKS.
Fig. \ref{2d-dambreak-mesh} shows the computational domain and the mesh.
The length of the dam breach is $75$ and it starts at $y=95$. The dam itself has a width of $10$ and its left side is located at $x = 95$.
At $t=0$ the stationary water surface has a discontinuity with $h_l=10$ and $h_r=\epsilon$ across the breach, and two values of $\epsilon=5$ and $\epsilon=1\times 10^{-3}$ are used to simulate the wet and dry bed cases, respectively.
For the wet case, the individual water levels of layer 1 and layer 2 are set as
\begin{equation*}
(h_1,h_2) = \begin{cases}
(9,1),  ~0\leq x<95,\\
(5,0),  ~95\leq x.
\end{cases}
\end{equation*}
For the dry case, the individual water levels of layer1 and layer 2 are set as
\begin{equation*}
(h_1,h_2) = \begin{cases}
(9,1),  ~0\leq x<95,\\
(\epsilon,\epsilon),  ~95\leq x.
\end{cases}
\end{equation*}
The boundary condition on the far right is the free boundary, and the other boundary conditions are the non-penetration slip wall boundaries.
The mesh size far from the breach is $h_{mesh}=2.5$, and is locally refined by $3.3$ times around the dam breach.

The 3-D water surface heights at $t=7.2$ are shown in Fig.\ref{2d-dambreak-wet} and Fig.\ref{2d-dambreak-dry}.
The discontinuous bore waves are captured without spurious oscillation.
It clearly shows that the wave propagating speed is higher in the dry bed case.

\section{Conclusion}

In this study, we have developed a compact high-order Gas-Kinetic Scheme (GKS) on a triangular mesh to solve the Two-Layer Shallow Water Equations (TLSWE). The compact scheme is highly accurate and robust in capturing discontinuous solutions.

The gas evolution model at the cell interface in the kinetic scheme explicitly captures the dynamics from the particle free transport, collisions, and acceleration from the external forcing term on the particle trajectory.
The time-accurate evolution solution provides not only the flow variable update inside each cell, but also the gradients of the flow variables.
As a result, based on the updated flow variables and their gradients, compact stencil can be used in the reconstruction and the design of the compact scheme.

The compact GKS has several key features in solving TLSWE.
The high-order compact reconstruction on a triangular mesh is naturally obtained. The existence of the time derivative of the flux function uses less stages to get a scheme with high-order accuracy in time, such as the two-stages for the fourth-order time accuracy.
This compact GKS provides accurate numerical solutions for the TLSWE and is ready for its engineering application in the coastal area ocean flow.

\section*{Acknowledgments}
The current research is supported by CORE as a joint research centre for ocean research between QNLM and HKUST through the project QNLM20SC01-A and QNLM20SC01-E, the National Natural Science Foundation of China (No. 12172316),
and Hong Kong research grant council 16208021 and 16301222.

\section*{References}
\bibliographystyle{ieeetr}
\bibliography{compact-gks}

\end{document}